# Non-negative diffusion bridge of the McKean–Vlasov type: analysis of singular diffusion and application to fish migration

Hidekazu Yoshioka[1, *]

[1] Graduate School of Advanced Science and Technology, Japan Advanced Institute of Science and Technology, 1-1 Asahidai, Nomi, Ishikawa, Japan

* Corresponding author: yoshih@jaist.ac.jp, ORCID: 0000-0002-5293-3246

**Abstract**

The objective of this paper is to provide a new mathematical tool for fish migration that has not been studied well. McKean–Vlasov stochastic differential equations (MVSDEs) have broad potential applications in science and engineering, but remain insufficiently explored. We consider a non-negative McKean–Vlasov diffusion bridge, a diffusion process pinned at both initial and terminal times, motivated by diurnal fish migration phenomena. This type of MVSDEs has not been previously studied. Our particular focus is on a singular diffusion coefficient that blows up at the terminal time, which plays a role in applications of the proposed MVSDE to real fish migration data. We prove that the well-posedness of the MVSDE depends critically on the strength of the singularity in the diffusion coefficient. We present a sufficient condition under which the MVSDE admits a unique strong solution that is continuous and non-negative. We also apply the MVSDE to the latest fine fish count data with a 10-min time interval collected from 2023 to 2025 and computationally investigate these models. Thus, this study contributes to the formulation of a new non-negative diffusion bridge along with an application study.

***Statements & Declarations***

**Acknowledgments:** The author thanks the Japan Water Agency, Ibi River and Nagara River General Management Office, for providing valuable fish migration data in recent years.

**Funding:** This study was supported by the Japan Science and Technology Agency (PRESTO No. JPMJPR24KE).

**Conflict of Interests:** The author declares no conflicts of interest.

**Data Availability:** Data will be made available upon reasonable request to the corresponding author.

**Declaration of Generative AI in Scientific Writing:** The author did not use generative AI in the scientific writing of this manuscript.

**Contribution:** The author prepared all parts of this manuscript.

## 1. Introduction

### 1.1 Research background

Fish migration is a complex biological phenomenon closely linked to human life. Numerous migratory fish species, such as salmonids, have high commercial value, and assessing their biomass is crucial [1]. Biomass assessment is also important for the conservation of endemic migratory fish species [2]. Physical barriers, such as weirs and dams, critically affect the life cycles of migratory fish species, and studying their migration behavior is of pivotal importance for sustainable river environmental planning [3, 4]. However, the dynamics of migratory fish species, such as the determinacy and randomness of dispersal, are largely unknown [5]. The migration dynamics of fish as fishery resources affects market-making [6]. Thus, fish migration has become an important research topic in the modern world.

Modeling migration of animals, including fish species, have been addressed by using stochastic differential equations (SDEs) [7] and related mathematical models because they can phenomenologically describe complex phenomena [8]. Stochastic agent-based models have been applied to the computational study of transient space-use dynamics [9]. SDEs with boundary reflections have been applied to animal movements in confined domains [10]. SDEs that govern migrating fish populations have been proposed in the literature, where the speed or timing of migration is assumed to arise from an optimality principle [11, 12], and daily fish counts have been studied using jump-driven SDEs [13, 14]. Optimization formalism has also been applied to bird migration [15]. Fish migration can be considered as population dynamics, where the decisions of individuals arise as those of the population. From this perspective, mean-field games based on diffusion processes have been effectively employed to determine the optimal (ideal-free) distribution considering environmental capacity [16] and social learning [17]. Mean-field games have also been applied to other collective phenomena, such as decision-making in honeybees [18] and swarm movements with collision avoidance [19].

From a mathematical viewpoint, SDEs whose coefficients depend on the mean-field (e.g., probability density or moments of the variables governed by SDEs), as found in mean-field games, are called McKean–Vlasov SDEs (MVSDEs). The fundamental theory of MVSDEs, such as the existence and uniqueness [20], long-time behavior of solutions [21], cases that are fractional in time [22], and cases with jumps [23, 24] have been studied well, while their applications to real problems remain limited except for a few study areas, such as finance and economics [25-27] and machine learning [28-30].

Recent technological progress has enabled the observation of fish migration at temporally fine scales (i.e., intraday scales, such as hourly or even sub-hourly scales), significantly advancing our understanding of the biology and ecology of migratory fish. The hourly migration of eels has been studied with a focus on the safe management of tidal barrages [31] and navigation locks [32]. The relationship between the tidal and circadian rhythms of the European eel has been quantified in narrow seas [33]. The passage probability of Atlantic salmon along a heavily human-influenced river has been estimated at a sub-daily resolution [34]. The intraday bidirectional migration of smolting salmon has been studied using acoustic biologging, suggesting that this phenomenon is related to the presence of mammalian predators [35]. The spawning migration of rheophilic fish has been found to have daily periodicities that vary according to sex [36]. The thermal habitat selection and vertical activity of walleyes in a river-lake connected system have been clarified at both seasonal and intraday scales [37]. The performance of intraday fishways in large anthropized rivers has been evaluated for more than ten fish species [38]. Hourly environmental DNA sampling was conducted to clarify the intraday spawning activities of *Plecoglossus altivelis altivelis* (Ayu, sweetfish) in Japan [39], which is a major inland fishery resource in the country.

Despite the progress in studies on temporally-fine fish migration phenomena, mathematical modeling has not been well studied. Recently, Yoshioka [40] proposed the application of a non-negative diffusion bridge, which is an SDE with a unique non-negative solution and pinned initial and terminal conditions, to modeling the intraday fish count of *Ayu* between sunrise and sunset during the upstream migration season. However, this SDE does not have any mean-field components despite describing an aspect of fish population dynamics; moreover, it has been proposed heuristically, and how individual behavior emerges as the SDE remains unclear. Diffusion bridges without mean-field components have been studied by focusing on optimization principles [41, 42] and conditioning [43, 44]. To the best of the author's knowledge, bridge-type MVSDEs have only been studied for Brownian bridge cases [45]. Mathematically studying non-negative diffusion bridges that have mean-field components and exploring their individual-based origin would deepen our understanding of fish migration phenomena and enable us to develop a more realistic theory for describing them.

### 1.2 Aim and contribution

The aims of this study are to propose an MVSDE motivated by the non-negative diffusion bridge [40] for intraday fish migration and present its application to the latest real dataset. Therefore, this study analyzes a novel mathematical model and its application to deepen both fundamental and applied research areas focusing on fish migration. The following section outlines our contributions to the aims of this study.

The MVSDE considered in this study is of the bridge type, with the process pinned at both the initial and terminal conditions. This is a nonlinear SDE having mean-field components. Its drift and diffusion coefficients originate from the classical Cox–Ingersoll–Ross (CIR) process of interest rate dynamics in finance [46]. The classical CIR process admits a global unique non-negative strong solution, with explicitly known moments and characteristic functions [47]. The analytical tractability of the CIR process has been utilized not only in economics [48] but also in other research areas, such as insurance [49], chemistry [50], cell biology [51, 52], fluid dynamics [53], and hydrological processes [54]. The proposed MVSDE is reduced to the classical CIR process if the pinning at the terminal time is abandoned and the coefficients are simplified. We show that our MVSDE partly inherits the analytical tractability of the classical CIR process. Moreover, our MVSDE is a generalization of the non-negative diffusion bridge [40] because the latter is obtained by neglecting the mean-field effect and simplifying the drift and diffusion coefficients. We show that the MVSDE admits a unique strong solution pinned at the initial and terminal times with probability 1, along with a sufficient condition for the singularity of the diffusion coefficient, which plays a key role in our application study.

From a biological standpoint, the mean-field nature of the proposed MVSDE arises from its derivation from an individual-based model, because fish migration is a biological phenomenon in which many individuals move simultaneously. We present a bottom-up derivation procedure of the MVSDE by relying on the regeneration of CIR processes (e.g., Exercise 1.2.13 in [47]), i.e., adding up certain CIR processes results in another CIR process. Despite its phenomenological model, the proposed MVSDE has an individual-based origin. We also note that diffusion bridges, such as Brownian bridges, have been applied to the spatial dynamics of fish (e.g., [55-57]), whereas our focus is on fish count time series data at a fixed point; therefore, the modeling objectives are different from each other.

Finally, we apply the proposed MVSDE to the latest dataset of juvenile upstream migration of *Ayu* from 2023 to 2025, in which migrants are counted in 10-min intervals between sunrise and sunset each day in spring to early summer. The temporally fine, intraday upstream migration behavior of this fish has only recently begun to be studied [13, 58]. The initial and terminal conditions of the MVSDE corresponded to sunrise and sunset, respectively, during which the migratory behavior of the fish is activated. This dataset suggests that intraday fish migration behavior is intermittent. The analytical tractability of the model, in which the time-dependent average and variance are explicitly available, is fully exploited in the model identification procedure of this study. The identified models fit the data well and can generate intermittent time series according to the application of the Feller condition, along with numerical simulation. We also computationally investigate the case where the singularity of the diffusion coefficient is sufficiently strong to violate the sufficient condition for the unique existence of solutions to our MVSDE.

### 1.3 Structure of this paper

The remainder of this paper is organized as follows. **Section 2** presents the MVSDE and examines its well-posedness. Its individual-based origin is explained and a sufficient condition for the unique existence of strong solutions to the MVSDE, which is the main theoretical result of this study, is presented. **Section 3** applies the MVSDE to unique 10-min fish migration data. A computational study of the MVSDE is also presented in this section. **Section 4** summarizes the study and discusses topics to be addressed in future. The **Appendix** presents the proofs and auxiliary data.

## 2. McKean–Vlasov SDE of the bridge type

### 2.1 Formulation

We use a complete probability space, as standard in stochastic calculus [59]. Throughout this study, time is a continuous variable denoted by $t$. The expectation is denoted by $\mathbb{E}$. We consider the following MVSDE that governs the stochastic process $X=\left(X_t\right)_{0\leq t\leq T}$ in a time interval $\left[0,T\right]$ with a fixed terminal time $T>0$:

$$\underset{\text{Increment}}{\mathrm{d}X_t} = \underbrace{\left(\underbrace{a\left(t,\bar{X}_t\right)}_{\text{Source}} - \underbrace{\frac{r\left(t,\bar{X}_t\right)}{T-t}X_t}_{\text{Mean reversion}}\right)\mathrm{d}t}_{\text{Drift}} + \underbrace{\sigma\left(t,\bar{X}_t\right)\sqrt{\frac{r\left(t,\bar{X}_t\right)}{\left(T-t\right)^{\alpha}}X_t}\,\mathrm{d}B_t}_{\text{Diffusion}},\quad 0<t<T \tag{1}$$

subject to the pinned initial condition $X_0 = 0$ and terminal condition $X_T = 0$. Here, $\bar{X} = \left(\bar{X}_t\right)_{0\le t\le T}$ denotes the expectation of $X$, that is, $\bar{X}_t = \mathbb{E}\left[X_t\right]$; $a:\left[0,T\right]\times\mathbb{R}\to\left[0,+\infty\right)$ denotes the source coefficient that is continuous; $r:\left[0,T\right]\times\mathbb{R}\to\left[0,+\infty\right)$ indicates the reversion coefficient that is continuous; $\sigma:\left[0,T\right]\times\mathbb{R}\to\left[0,+\infty\right)$ denotes the volatility that is continuous and controls the noise intensity; $\alpha\in\mathbb{R}$ represents a singularity parameter of diffusion; and $B = \left(B_t\right)_{0\le t\le T}$ represents a standard 1-D Brownian motion. The diffusion term in (1) is defined in the Itô's sense [7]. In the rest of **Section 2**, we set $T = 1$ without any loss of generality:

$$\mathrm{d}X_t = \left(a\left(t,\bar{X}_t\right) - \frac{r\left(t,\bar{X}_t\right)}{1-t}X_t\right)\mathrm{d}t + \sigma\left(t,\bar{X}_t\right)\sqrt{\frac{r\left(t,\bar{X}_t\right)}{\left(1-t\right)^{\alpha}}X_t}\,\mathrm{d}B_t,\quad 0<t<1. \tag{2}$$

The MVSDE (1) reduces to the non-negative diffusion bridge of Yoshioka [40] if the coefficients $a$, $r$, and $\sigma$ are set to positive constants and $\alpha$ is set to 1, thus representing a generalized version of the previous model with a generalized singularity in the diffusion coefficient. From the viewpoint of fish migration, the variable $X$ represents unit-time fish count (e.g., number of migrants in each fixed time interval) observed at a fixed point along a river. In the context of diurnal migration, the initial and terminal times correspond to sunrise and sunset, respectively. The MVSDE (1) then describes temporal evolution of the unit-time fish count by accounting for its deterministic (drift) and stochastic (diffusion) parts, both being singular at the terminal time. The mean-field component, denoted by $\bar{X}_t$ in the proposed model, represents an individual fish's prediction of fish count dynamics, as the expectation is interpreted as a coarse estimate of a stochastic process. The expectation in the drift and diffusion coefficients represents interactions among individual fish that eventually shape the fish migration phenomenon. It may be more realistic to use an expectation conditioned on the current state instead of the usual expectation $\bar{X}_t$; however, this would result in a far more complex model that potentially losses the analytical tractability of the proposed model; therefore, we do not adopt this method in this study.

Phenomenologically, the singularity of the drift (the second term in the drift term of (1)) enforces the average $\bar{X}_t = \mathbb{E}\left[X_t\right]$ to vanish at the terminal time. The singularity of the diffusion then allows for modeling temporal dependence of variance $\mathbb{V}\left[X_t\right] = \mathbb{E}\left[\left(X_t - \bar{X}_t\right)^2\right]$. The case $\alpha = 1$ originates from the introduction of a biological clock, time change in stochastic analysis, into the model [40]. We generalize this singularity to the case $\alpha \ne 1$, and the cases $\alpha > 0$ and $\alpha < 0$ amplifies and dampens the diffusion at the terminal time, respectively. The possible range of $\alpha$ needs to be determined considering other coefficients in the MVSDE (1). More specifically, we show that $\alpha$ must be sufficiently small, such that the blow-up of the diffusion near the terminal time is suppressed by the reversion mechanism in the drift. Otherwise, the pinning at the terminal time cannot be prescribed continuously in time (see **Section 3** and **Appendix C** for an explicit model), and the well-posedness of the MVSDE may be broken. Parameters in the MVSDE can be determined by fitting it to time series data.

### 2.2 Unique existence and moments

We present a unique existence result of strong solutions (pathwise solutions that are continuous on $[0,1]$ with probability 1) to the MVSDE (1) under a certain condition. Two factors need to be considered here: the mean-field components and parameter $\alpha$. Intuitively, the dependence on the mean-field components should be sufficiently regular to ensure that the average $\bar{X}$ is uniquely determined. Moreover, the parameter $\alpha$ need to be small to ensure that the singularity of the diffusion coefficient is not too strong. **Assumption 1** is used to formalize these conditions:

***Assumption 1***

***(Lipschitz continuity)*** $\left|f(t,x_1)-f(t,x_2)\right| \le L_f\left|x_1-x_2\right|$ *for all* $t\in[0,1]$ *and* $x_1,x_2\in\mathbb{R}$ *(3)*

*with some global constant* $L_f>0$*, where* $f$ *represents* $a$*,* $r$*, and* $\sigma$*.*

***(Bound of reversion)*** $0<\underline{r}\le r(t,x)\le\bar{r}<+\infty$ *and* $0\le a(t,x)\le\bar{a}<+\infty$ *for all* $t\in[0,1]$ *and* $x\in\mathbb{R}$ *(4)*

*with some global constants* $\underline{r},\bar{r},\bar{a}>0$.

***(Small singularity)*** $\alpha<\min\{2,1+\underline{r}\}$. *(5)*

The motivation for each item in **Assumption 1** is as follows. The Lipschitz continuity (3) controls the regularity of the coefficients. The bound of the reversion (3) controls sizes of the coefficients $a$ and $r$. Combining (3) and (4) yields that the following ordinary differential equation (ODE) admits a unique continuous solution $u=(u_t)_{0\le t\le 1}$:

$$\frac{\mathrm{d}u_t}{\mathrm{d}t}=a(t,u_t)-\frac{r(t,u_t)}{1-t}u_t,\quad 0<t<1 \tag{6}$$

with initial and terminal conditions of $u_0=u_1=0$. We expect that the expectation $\bar{X}$ is the unique solution to (6) and hence $u=\bar{X}$, which turns out to be true under **Assumption 1**. Then, one must be careful about the regularity of the diffusion term, particularly the parameter $\alpha$, to ensure that it induces a martingale continuous at the terminal time. Otherwise, the variance of $X$ does not vanish at the terminal time and the terminal condition is not satisfied with probability 1. The small singularity condition (5) prevents this issue.

Now, we state our main theoretical result.

***Proposition 1***

*Under **Assumption 1**, there exists a unique strong solution to the MVSDE (1) that satisfies the initial and terminal conditions* $X_0 = X_1 = 0$ *with probability 1.*

The small singularity condition (5) is sharp as shown in (36) in **Appendix** because it shows that the variance of $X$ does not vanish at the terminal time when $\alpha \geq 1+r$ or $\alpha \geq 2$. Without the small singularity condition, the existence of strong solutions is not lost in any $(0,1-\varepsilon)$ with small sufficiently $\varepsilon > 0$ but the left limit of the variance at the terminal time is not necessarily 0.

***Remark 1*** Inspecting **Proof of Proposition 1** suggests that **Assumption 1** can be weakened as follows: the coefficient $a$ is allowed to have a singularity as follows:

$$a(t,x) = a_0(t,x) + \frac{A}{1-t}, \text{ for all } t \in [0,1) \text{ and } x \in \mathbb{R} \tag{7}$$

with a real constant $A$ such that $A < \underline{r}$, where $a_0$ denotes a bounded and measurable function that satisfies (3). In this case, we need to replace the upper-bound of $\alpha$ in (5) by $\min\{2, 1+\underline{r}-A\}$. This implies that the source rate, biologically indicating that the supply of migrants from the downstream river of the observation point, is allowed to rapidly increase near the terminal time. Moreover, the regularity of the coefficients $a, r$ can be weakened, such as only locally Lipschitz continuous with respect to the second argument, provided that the ODE (6) admits a unique solution.

### 2.3 Individual-based origin

We begin with the following findings regarding CIR processes. Let $a_1, a_2, a_3 > 0$ be parameters, and $B^{(1)}$ and $B^{(2)}$ be independent standard 1-D Brownian motions. Let $X^{(1)}$ and $X^{(2)}$ be CIR processes governed by the following Itô's SDE

$$\mathrm{d}X_t^{(i)} = \left(\frac{a_1}{2} - a_2 X_t^{(i)}\right)\mathrm{d}t + a_3\sqrt{X_t^{(i)}}\mathrm{d}B_t^{(i)}, \quad t > 0, \quad i = 1,2 \tag{8}$$

subject to the initial conditions $X_0^{(i)} = \frac{x}{2} \geq 0$ ($i = 1,2$). Then, the summation $\hat{X} = X^{(1)} + X^{(2)}$ is a CIR process governed by the following SDE:

$$\mathrm{d}\hat{X}_t = \left(a_1 - a_2\hat{X}_t\right)\mathrm{d}t + a_3\sqrt{\hat{X}_t}\mathrm{d}\hat{B}_t, \quad t > 0, \quad i = 1,2 \tag{9}$$

subject to the initial condition $\hat{X}_0 = x$, where $B$ represents some standard 1-D Brownian motion. Indeed, direct calculations show

$$
\begin{aligned}
\mathrm{d}\hat{X}_t &= \mathrm{d}X_t^{(1)} + \mathrm{d}X_t^{(2)} \\
&= \left(\frac{a_1}{2} - a_2 X_t^{(1)}\right)\mathrm{d}t + a_3\sqrt{X_t^{(1)}}\mathrm{d}B_t^{(1)} + \left(\frac{a_1}{2} - a_2 X_t^{(2)}\right)\mathrm{d}t + a_3\sqrt{X_t^{(2)}}\mathrm{d}B_t^{(2)} \\
&= \left(a_1 - a_2\hat{X}_t\right)\mathrm{d}t + a_3\left(\sqrt{X_t^{(1)}}\mathrm{d}B_t^{(1)} + \sqrt{X_t^{(2)}}\mathrm{d}B_t^{(2)}\right) \\
&= \left(a_1 - a_2\hat{X}_t\right)\mathrm{d}t + a_3\sqrt{\hat{X}_t}\mathrm{d}\hat{B}_t
\end{aligned}
, \quad t > 0. \qquad (10)
$$

Here, the last equality is understood in the sense of law (e.g., Exercise 1.2.13 in [47]).

As demonstrated above, adding up CIR processes (8) with scaled source coefficients ($a_1$ divided by 2 above) yields another CIR process (9) with the source coefficient $a_1$. As stated in **Proposition 2** below, similar calculations apply to the MVSDE (1) by considering it a CIR process with time-dependent coefficients.

***Proposition 2***

*Under **Assumption 1**, for any* $n \in \mathbb{N}$, *consider the following SDEs*

$$
\mathrm{d}X_t^{(i)} = \left(\frac{a\left(t,\bar{X}_t\right)}{n} - \frac{r\left(t,\bar{X}_t\right)}{1-t}X_t^{(i)}\right)\mathrm{d}t + \sigma\left(t,\bar{X}_t\right)\sqrt{\frac{r\left(t,\bar{X}_t\right)}{\left(1-t\right)^{\alpha}}X_t^{(i)}}\mathrm{d}B_t^{(i)}, \quad 0 < t < 1, \quad i = 1,2,\dots,n \qquad (11)
$$

*subject to* $X_0^{(i)} = X_1^{(i)} = 0$ $(i = 1,2,\dots,n)$, *where* $B^{(i)}$ $(i = 1,2,\dots,n)$ *represent independent standard 1-D Brownian motions. Then, it follows that*

$$
X_t = \sum_{i=1}^{n} X_t^{(i)}, \quad 0 \le t \le 1 \qquad (12)
$$

*in the sense of law, where* $X$ *denotes the unique strong solution to the MVSDE (1) and* $\bar{X}$ *its expectation.*

**Proposition 2** suggests that the MVSDE (1) as a mathematical model of the unit-time fish count can be considered an aggregation of small, similar processes (11). Here, the word "small" arises because the source $a$, controlling the size of both average and variance of the process is divided by $n$ in (11). SDEs (11) that govern small processes can be considered as the governing equations of individual migration activities if we consider $n$ as the total number of individuals in the downstream area of the observation point for counting fish. From this perspective, the representation (12) allows the unit-time fish count to be viewed as an aggregated quantity derived from the individual-based model on the right-hand side. This individual-based formulation also applies to the non-negative diffusion bridge of Yoshioka [40] because our MVSDE is its generalization. Thus, this study provides a new theoretical characterization of diffusion bridges for describing fish migration phenomena.

***Remark 2*** Taking a macroscopic limit of an individual-based model that tracks swimming behavior of individual fish in parallel [60-63] would be another strategy to obtain an SDE of fish count at a fixed observation point; however, the extent of allowable complexity in such individual-based models is unclear.

The CIR process can also be obtained as a high-frequency limit of self-exciting jump processes [64, 65], suggesting another route to obtain the MVSDE, where each jump is considered a migration event of one individual fish with an adequate scaling argument.

***Remark 3*** We assume a homogeneous population in **Proposition 2** such that the source coefficient is the same for all $i = 1,2,...,n$. By inspecting **Proof of Proposition 2**, this assumption can be relaxed, and the conclusion of this proposition holds if SDEs (11) are replaced by

$$\mathrm{d}X_t^{(i)} = \left( a_i\left(t,\bar{X}_t\right) - \frac{r\left(t,\bar{X}_t\right)}{1-t} X_t^{(i)} \right)\mathrm{d}t + \sigma\left(t,\bar{X}_t\right)\sqrt{\frac{r\left(t,\bar{X}_t\right)}{\left(1-t\right)^{\alpha}} X_t^{(i)}}\,\mathrm{d}B_t^{(i)},\quad 0<t<1,\quad i=1,2,...,n \tag{13}$$

without changing initial or terminal conditions. Here, $a_i : [0,1]\times\mathbb{R} \to [0,+\infty)$ is a coefficient that satisfies the assumption of $a$ stated in **Proposition 1**, such that $a(t,x) = \sum_{i=1}^{n} a_i(t,x)$ for all $0 \le t \le 1$ and $x \in \mathbb{R}$. Therefore, the formula (12) also applies to this heterogeneous population case.

## 3. Application

### 3.1 Study target

The target fish species in this study is the amphidromous *Ayu* that migrate along the Nagara River, a class A river in the Kiso River system in the Tokai region of Japan (**Figure 1**). No large dams have been constructed along the mainstem that would fully obstruct fish migration. Pristine water quality of the Nagara River, as well as the surrounding cultural and ecological significance related to *Ayu* fisheries, led to its designation as a Globally Important Agricultural Heritage System by the United Nations Educational, Scientific and Cultural Organization in 2015 [66].

The amphidromous *Ayu* migrates between a river and the sea (coastal area) and has a one-year life history. Hence, they do not have an age structure. The life history of *Ayu* can be summarized as follows [67]: spawning of mature fish occurs in rivers in autumn, after which the next generation begins. After approximately two weeks of spawning, the hatched larvae, approximately 6 mm long, drift downstream to the sea connected to the river. These larvae spend winter in coastal areas of the sea. In the spring, juveniles with lengths of approximately 60 mm begin their upstream migration to nearby rivers. They grow to approximately 150–300 mm by the end of summer in the river. Mature fish spawn and die. In the Nagara River, the juvenile upstream migration of fish is observed from March to June each year.

The 10-min count data of juvenile *Ayu* have been acquired automatically from a fishway installed in the Nagara River Estuary Barrage since 2021[1], and the dataset from 2023 to 2025 was available. To the best of the author's knowledge, this type of temporally-fine intraday data on fish migration is not always available (at least for research purposes). The dataset used here adds data for 2025 to those for 2023 and 2024 studied in [40]. The pairs (start, end) of observation dates are (February 22, June 30) in 2023, (February 26, June 30) in 2024, and (March 2, June 30) in 2025. In total, 852,596, 1,236,102, and 1,318,281 juveniles were counted in 2023, 2024, and 2025, respectively. **Figure 2** shows the 10-min fish count data for each year. The 10-min dataset suggests that the upstream migration of *Ayu* is an intermittent phenomenon with burst-and-rest structures, as highlighted in [40] for 2023 and 2024. **Figure 2** implies that this finding also applies to the data in 2025. This point will be theoretically studied in terms of a Feller condition later.

[1] Japan Water Agency, Ibi River and Nagara River General Management Office. https://www.water.go.jp/chubu/nagara/15_sojou/r07_ayu_suuchi.html (last accessed on September 21, 2025)

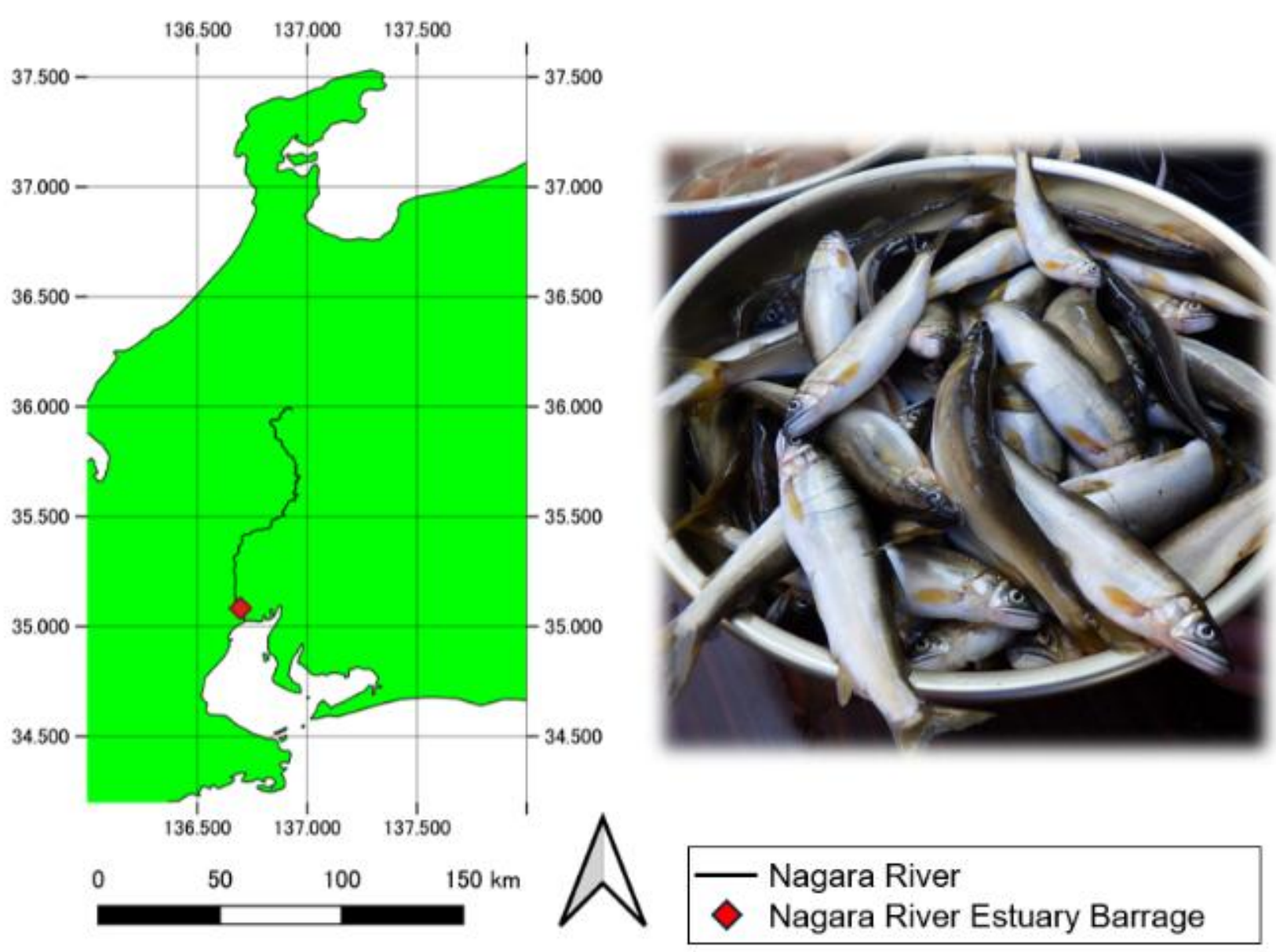

**Figure 1.** Map of the study site (left) and the target fish species, *Ayu* (right).

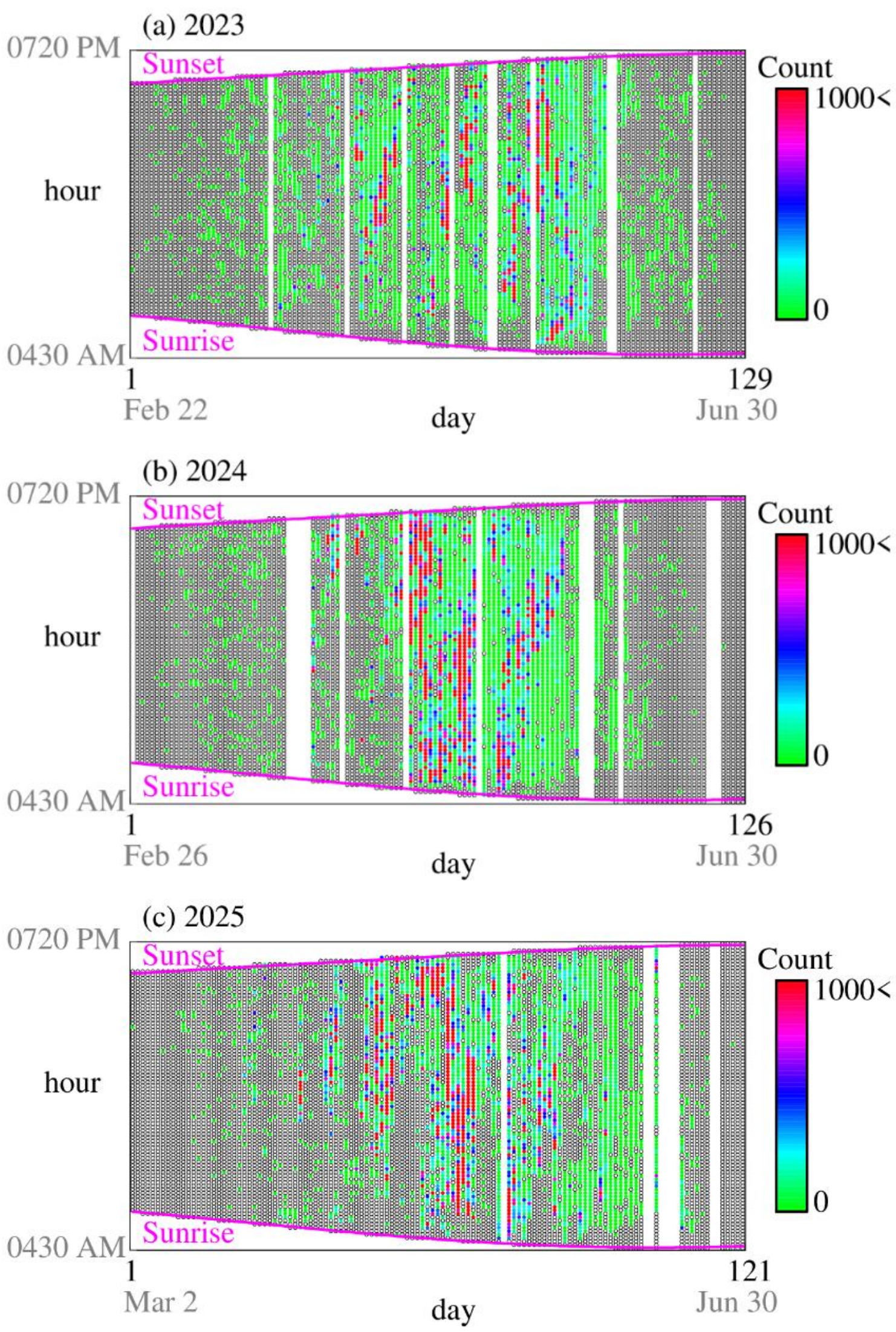

**Figure 2.** 10-min fish count of *Ayu* in (a) 2023, (b) 2024, and (c) 2025. White circles represent zero count data, and blank areas represent no data (The figure panels (a) and (b) were adapted from Yoshioka [40] with modifications).

### 3.2 Data normalization and analysis

The normalization procedure follows that of Yoshioka [40] with a proper modification to address the singularity in the diffusion during $\alpha \neq 1$. In this subsection, the arguments $\left(t, \bar{X}_t\right)$ are omitted from coefficients to simplify the notation.

Sunrise is set at the initial time 0 each day. Assume that there are $M \in \mathbb{N}$ observation days in a year (e.g., 2023, 2024, 2025). The sunset on Day $N$ ($1 \leq N \leq M$) in the year is denoted as $T_N$. The total number of migrants on Day $N$ is denoted as $Y_N$, which is assumed to be positive here. The time increment (10 min) of the data is denoted as $\Delta t$. The MVSDE describing intraday fish migration on Day $N$ is expressed as follows:

$$\mathrm{d}X_{t,N} = \left(a_N - \frac{r_N}{T_N - t} X_{t,N}\right)\mathrm{d}t + \sigma_N \left(T_N - t\right)^{\frac{1-\alpha}{2}} \sqrt{\frac{r_N}{T_N - t} X_{t,N}}\,\mathrm{d}B_{t,N}, \quad 0 < t < T_N \tag{14}$$

with the initial condition $X_{0,N} = 0$ and terminal condition $X_{T_N,N} = 0$, where the subscript $N$ represents the quantities on Day $N$, and each $B_{t,N}$ is assumed to be independent for different $N$. Because the total fish count $Y_N$ is different for different days, we employ the following normalized variables, all being non-dimensional: $s_N = \frac{t}{T_N} \in [0,1]$, $Z_{s,N} = \frac{\Delta t}{Y_N} X_{t,N}$, $\bar{a} = \frac{a_N T_N \Delta t}{Y_N}$, and $\bar{\sigma} = \sigma_N T_N^{\frac{\alpha-1}{2}} \sqrt{\frac{\Delta t}{Y_N}}$ by assuming the existence of common coefficients $\bar{a}, \bar{\sigma} > 0$ and $\bar{r} = r_N$ for all $N$. Substituting these quantities into (14) yields the following normalized MVSDE:

$$\begin{aligned} \mathrm{d}\left(\frac{Y_N}{\Delta t} Z_{s,N}\right) &= \left(\bar{a}\frac{Y_N}{T_N \Delta t} - \frac{\bar{r}}{T_N - T_N s}\frac{Y_N}{\Delta t} Z_{s,N}\right)\mathrm{d}\left(T_N s\right) \\ &\quad + \left(\bar{\sigma} T_N^{\frac{\alpha-1}{2}} \sqrt{\frac{Y_N}{\Delta t}}\right)\left(T_N - T_N s\right)^{\frac{1-\alpha}{2}} \sqrt{\frac{\bar{r}}{T_N - T_N s}\frac{Y_N}{\Delta t} Z_{s,N}}\,\mathrm{d}\left(\sqrt{T_N} B_{s,N}\right) \end{aligned}, \quad 0 < s < 1 \tag{15}$$

with $Z_{0,N} = Z_{s,N} = 0$. The term $\sqrt{T_N} B_{s,N}$ in (15) is understood under the transformation $B_{t,N} \to \sqrt{T_N} B_{s,N}$ in the sense of law with some abuse of notations. We thus obtain

$$\mathrm{d}Z_{s,N} = \left(\bar{a} - \frac{\bar{r}}{1-s} Z_{s,N}\right)\mathrm{d}s + \bar{\sigma}\sqrt{\frac{\bar{r}}{(1-s)^{\alpha}} Z_{s,N}}\,\mathrm{d}B_{s,N}, \quad 0 < s < 1. \tag{16}$$

This is the MVSDE (1) with terminal time $T = 1$. Based on the normalization procedure described above, we directly fit the terminal time $T = 1$ to the normalized intraday fish count data for each year. With the normalization procedure explained above, each intraday of the time-series data with $Y_N > 0$ can be considered one sample path of the MVSDE. We do not use the data from days on which no migrants are observed (i.e., $Y_N = 0$).

#### 3.2.1 Identified models

We fit MVSDE (1) for each year (2023, 2024, and 2025 separately) as well as for the aggregated years (2023–2025 in a aggregated way). The last case is examined because it can increase the total number of empirical sample paths under the normalization hypothesis for intraday fish migration, as explained above.

We consider the following model specifications with which the average and variance are explicitly found, as shown in **Appendix C**:

$$a\left(t,\bar{X}_t\right)=a>0 \quad \text{and} \quad r\left(t,\bar{X}_t\right)=r>0 \quad \text{(positive constants)} \tag{17}$$

and

$$\sigma\left(t,\bar{X}_t\right)=\sqrt{\mu^2+\omega\bar{X}_t}>0 \quad \text{(expectation-dependent coefficient)} \tag{18}$$

with parameters $\mu>0$ and $\omega\in\mathbb{R}$ chosen such that the right-hand side of (18) is well-defined as a positive value. A motivation behind this model specification arises from the earlier result [40] that the diffusion bridge without mean-field components does not fit the observed standard deviation satisfactorily. We considered that employing a more sophisticated diffusion term would improve this issue. The sign of the parameter $\omega$ determines how the mean-field effect affects the fish migration; the positive $\omega$ case implies that the mean-field effect amplifies the stochasticity involved in the fish migration, and vice versa.

For the MVSDE (1) in this case, parameters to be estimated are $a,r,\mu,\omega,\alpha$. The model identification procedure used in this study follows a two-step approach. First, the parameters $a,r$ are estimated via least squares by fitting the empirical and theoretical averages $\mathbb{E}\left[X_t\right]=\bar{X}_t$. Second, the remaining parameters $\mu,\omega,\alpha$ are estimated, also via least squares, by fitting the empirical and theoretical standard deviations $\sqrt{\mathbb{V}\left[X_t\right]}$. The advantage of this identification method is that it fully exploits the analytical tractability of the diffusion bridge, enabling explicit determination of the average and variance, and hence the standard deviation. No parameter constraints were not used for the parameter estimation in this paper.

### 3.3 Results and discussion

#### 3.3.1 Identified model

**Table 1** summarizes estimated parameter values for each case. We also examine the non-mean-field case, in which we specify $\omega=0$ *a priori* as shown in **Table 2**. Note that in the present setting, the mean-field component only affects stochasticity, and hence does not alter the average. In **Tables 1 and 2**, the normalized root mean-square error (RMSE) refers to the RMSE divided by the theoretical average integrated over the time interval $(0,1)$. All the identified values of $\alpha$ satisfy **Assumption 1**; hence, the corresponding MVSDEs are well posed according to **Proposition 1**. In all the cases, the condition (5) is satisfied. Also, the positivity of $\sigma^2$ is satisfied without any constraints in the simulation due to the parameter values in **Tables 1 and 2.**

**Figures 3 and 4** show the empirical and fitted averages and standard deviations, respectively, for each case. The theoretical results fit the empirical results reasonably well for each computational case.

As shown in **Tables 1 and 2,** incorporating the mean-field components improves the normalized RMSE, although the improvement is case-dependent; there is over a 10% reduction of the normalized RMSE for the models in 2023 and 2023–2025, whereas the reduction is a few percent for the models in 2024 and 2025. The normalized RMSEs for models with constant $a, r, \sigma$ with $\alpha = 1$, which are "Model 1" in Yoshioka [40] are presented in **Table B1** in **Appendix B**. A comparison between **Tables 2 and B1** shows that allowing $\alpha \neq 1$ reduces the normalized RMSE by approximately 10%, except for case 2023, where the reduction is less than 5%.

A key finding for the mean-field case is that all estimates yield $\omega < 0$, suggesting that the mean-field effect reduces stochasticity and hence randomness in fish migration. Owing to the CIR-like form of the proposed MVSDE, the decrease in diffusion reduced the intermittency of the fish count, as qualitatively discussed in the next subsection. Considering that intraday fish migration is described by a time-inhomogeneous SDE, it can be interpreted as a non-equilibrium phenomenon in which the interactions among individuals dynamically change. The theoretical finding that strengthening mean-field effects reduces the randomness in population dynamics has also been found in an individual-based model where the mean-field effects are correlated with the number of the topological neighborhood [68]. Moreover, the strength of polarization (i.e., strength of mean-field effects) was suggested to be increasing with respect to the average cluster size [69]; this qualitative relationship in our context may be understood as a less intermittent behavior of fish count where the total number of clusters is smaller and fish migrate collectively by forming a large wave. In these views, the identified models are considered consistent with the existing studies using different mathematical models. A recent study revealed that diurnal and seasonal social interactions exist in wild carp [70], whereas the existence of such a mechanism in *Ayu*, particularly during its upstream migration, remains an open issue. Addressing this issue from both theoretical and field standpoints is necessary to better quantify the mean-field effect on intraday fish migration behavior. A recent modeling study on collective fish schooling suggested a burst-and-rest pattern in an individual's swimming speed [71]. Stability of multi-agent system dynamics [72] may also affect the overall randomness in the fish migration. Analyzing how intermittency in the fish count changes depending on the species and environment would be a feasible way to address this issue.

Another important finding is that the estimated parameter values of $\alpha$ are approximately 0.3 to 0.5 and are therefore smaller than the fixed value 1 assumed in the previous model [40]. The models without the mean-field effect estimated smaller $\alpha$ values, thus supporting this consideration. Thus, the fitted models imply a weaker singularity of the diffusion term that was originally assumed [40]; nevertheless, all the mean-field cases suggest that there should be some degree of singularity in the diffusion term to better fit the empirical standard deviations. Regarding the theoretical standard deviations, **Figure 4** shows that the mean-field effect temporally predicts non-monotonic standard deviations in 2023 and 2023–2025. At this stage, it is difficult to determine whether this non-monotonicity is due to overfitting or is representing some inherent property of the real data due to fluctuations in the empirical standard deviations. Nevertheless, incorporating the mean-field effect opens the door towards a more flexible stochastic modeling of the juvenile upstream migration of *Ayu* at the study site.

**Table 1.** Estimated parameter values and normalized RMSEs for the mean-field case.

| Parameter or index | | 2023 | 2024 | 2025 | 2023-2025 |
|---|---|---|---|---|---|
| $a$ | | 6.003.E-02 | 3.219.E-02 | 3.021.E-02 | 3.673.E-02 |
| $r$ | | 1.647.E+00 | 5.137.E-01 | 4.574.E-01 | 7.100.E-01 |
| $\mu$ | | 1.500.E+00 | 1.483.E+00 | 1.613.E+00 | 1.634.E+00 |
| $\omega$ | | -1.229.E+02 | -9.898.E+01 | -1.298.E+02 | -1.439.E+02 |
| $\alpha$ | | 4.922.E-01 | 3.157.E-01 | 5.783.E-01 | 5.482.E-01 |
| Normalized RMSE | Ave | 3.818E-01 | 4.521E-01 | 5.092E-01 | 3.292E-01 |
| | Std | 4.438E-01 | 5.288E-01 | 7.613E-01 | 4.367E-01 |

**Table 2.** Estimated parameter values and normalized RMSEs for the non-mean field case ($\omega = 0$ is specified *a priori*).

| Parameter or index | | 2023 | 2024 | 2025 | 2023-2025 |
|---|---|---|---|---|---|
| $a$ | | 6.003.E-02 | 3.219.E-02 | 3.021.E-02 | 3.673.E-02 |
| $r$ | | 1.647.E+00 | 5.137.E-01 | 4.574.E-01 | 7.100.E-01 |
| $\mu$ | | 8.058.E-01 | 1.191.E+00 | 1.347.E+00 | 1.195.E+00 |
| $\alpha$ | | 2.176.E-01 | -1.284.E-01 | -9.916.E-02 | -1.218.E-01 |
| Normalized RMSE | Ave | 3.818E-01 | 4.521E-01 | 5.092E-01 | 3.292E-01 |
| | Std | 5.763E-01 | 5.555E-01 | 7.847E-01 | 4.972E-01 |

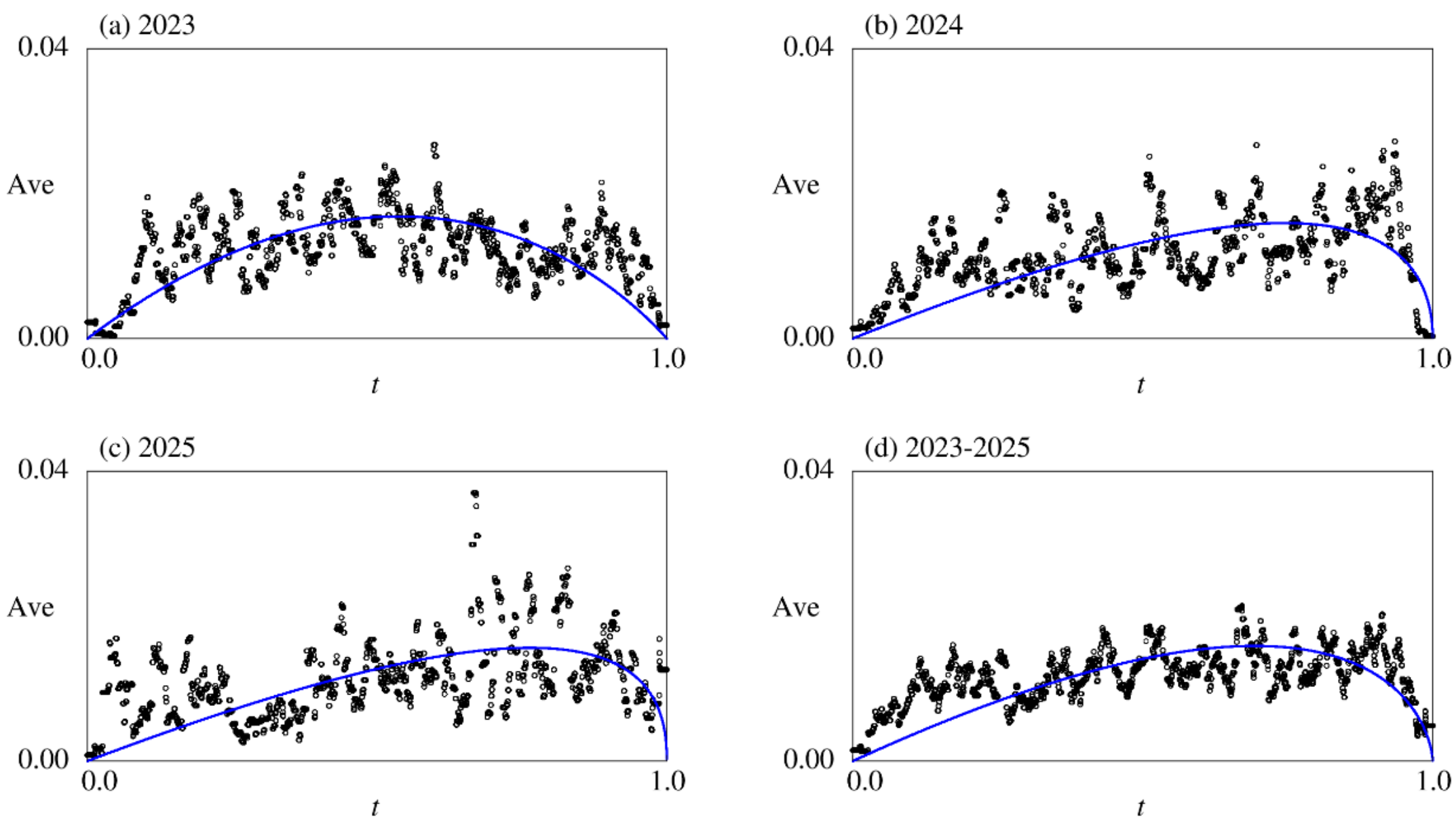

**Figure 3.** Comparison between empirical and fitted averages (Ave): (a) 2023, (b) 2024, (c) 2025, (d) 2023–2025. Circles and curves represent empirical and theoretical results, respectively. The averages are not affected by mean-field components.

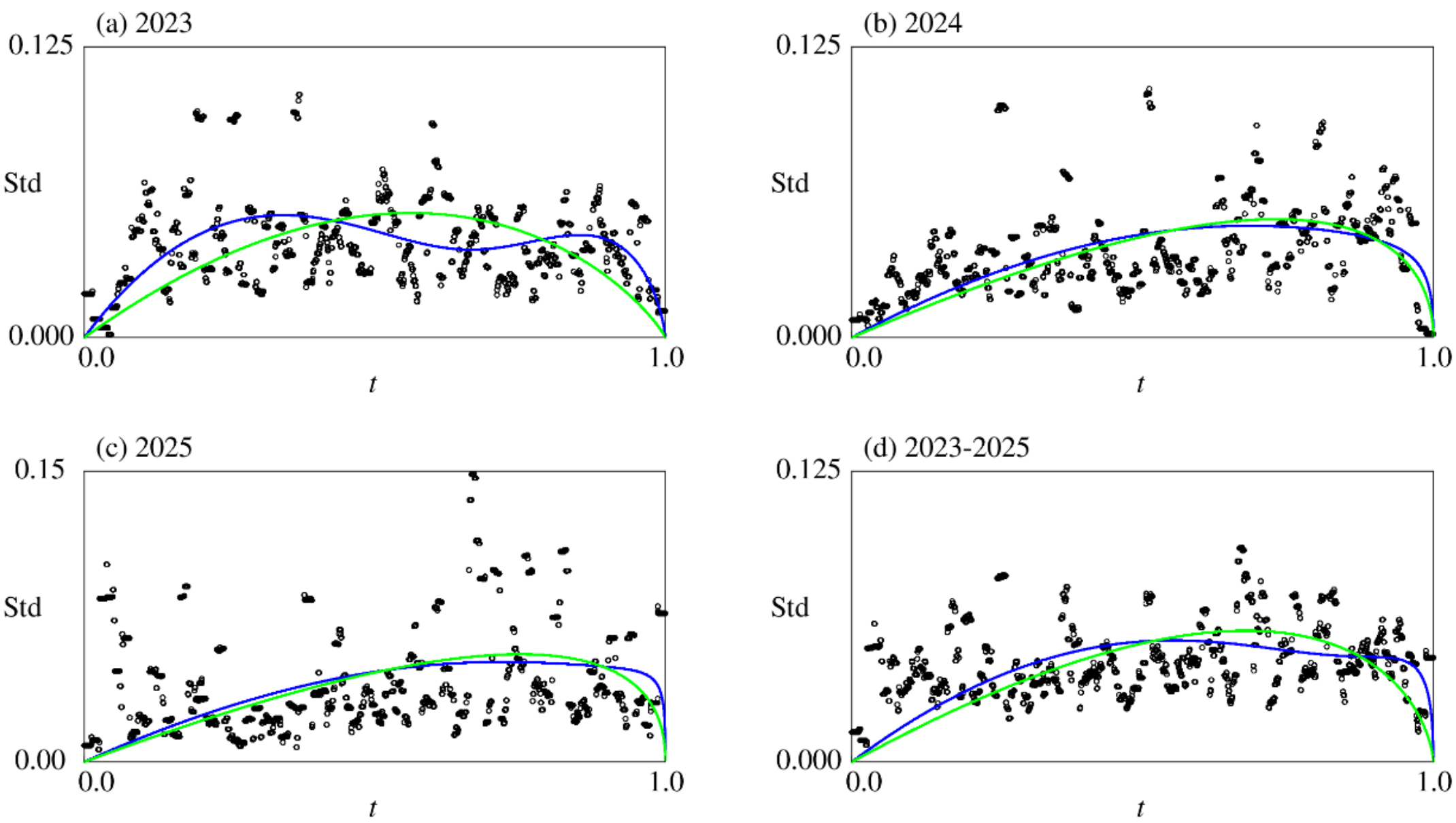

**Figure 4.** Comparison between empirical and fitted standard deviations (Std(a) 2023, (b) 2024, (c) 2025, (d) 2023–2025. Circles and curves represent empirical and theoretical results, respectively. Colors represent models with (blue) and without mean field effect (green). The standard deviations are affected by mean-field components.

### 3.3.2 Feller condition

The randomness and intermittency of the MVSDE (1) can be studied through the lens of the Feller condition, which is adapted here for the classical CIR process (e.g., Chapter 5.5.C in [73], Chapter 1.2.4 in [47]) with a modification of the time-dependent coefficients. To explain this condition, we introduce the non-dimensional index $F_t$, which measures the ratio between the noise intensity and acceleration in the drift (arguments of the coefficients are omitted to simplify the notation):

$$F_t = \frac{\left(\sigma\sqrt{\frac{r}{(1-t)^{\alpha}}}\right)^2}{a} - 1 = \frac{\sigma^2}{2a}\frac{r}{(1-t)^{\alpha}} - 1,\quad 0 < t < 1. \tag{19}$$

Intuitively, $F_t$ measures the dominance of diffusion $\sigma\sqrt{\frac{r}{(1-t)^{\alpha}}}$ relative to the source, $a$ at each time instance; a larger $F_t$ indicates higher dominance of diffusion, and vice versa. We say that the Feller condition is violated if $F_t \geq 0$ at all $0 \leq t \leq 1$, and is (partly) satisfied otherwise. Another important aspect of $F_t$ is that it, and hence the Feller condition, classifies regimes of the solution $X$ at each time instance; similar conditions have been discussed for various SDEs [74-78]. More specifically, in the present setting, the case $F_t \leq 0$ for some time interval $I = [t_1, t_2] \subset (0,1)$ corresponds to a low-volatility regime such that the solution $X_t$ ($t \in I$) never hits the time axis when $X_{t_1} > 0$, whereas the case $F_t > 0$ corresponds to a high-volatility regime under which $X_t$ ($t \in I$) possibly hits the time axis even when $X_{t_1} > 0$. Hence, the qualitative behavior of the solution to the MVSDE (1) can be classified by checking the sign of $F_t$. The Feller condition is also related to intermittency of sample paths of $X$ such that it has more intermittent sample paths under the high-volatility regime than under the low-volatility one. This is explained by the increased likelihood of $X$ hitting the time axis, i.e., the value 0, as fluctuations in the sample paths become larger.

According to the identified model, all cases violate the Feller condition, suggesting that they predict the migration of *Ayu* at the study site to be a random and intermittent phenomenon. **Figure 5** shows this finding computationally, where the average is theoretically the same in the two panels, whereas the intermittency is visually different owing to the distinctive Feller condition. This intermittent nature suggests the existence of clustering behavior behind fish migration. In the next subsection, we examine this point in more detail from a computational perspective.

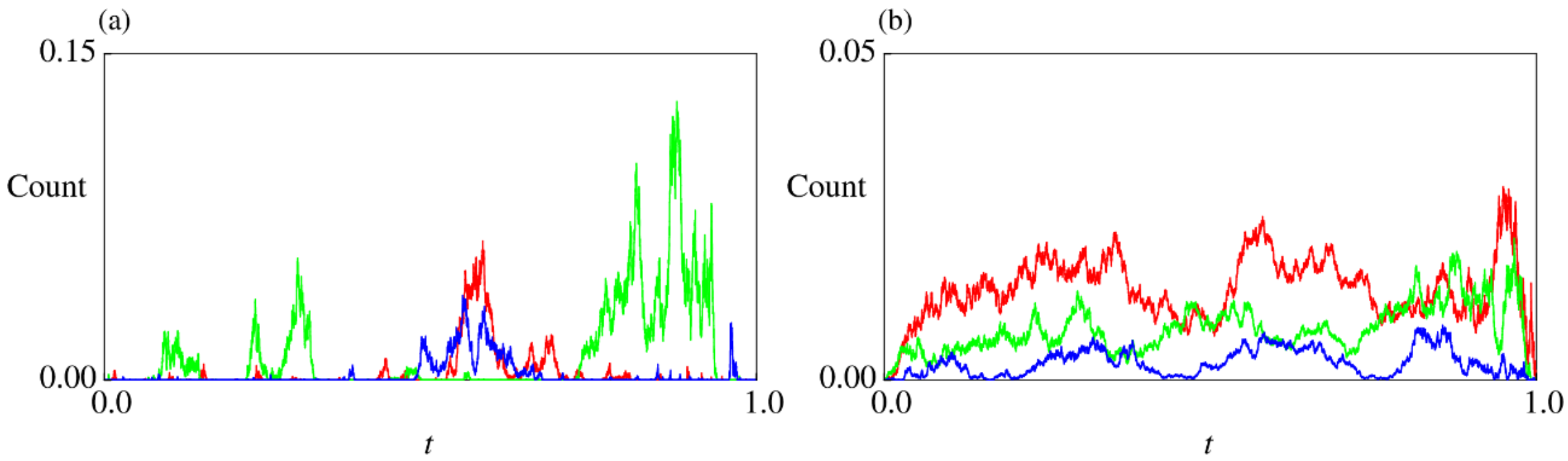

**Figure 5.** Computed sample paths of the unit time fish count: (a) identified model for 2023–2025 and (b) identified model for 2023–2025 where an artificial factor 0.2 is multiplied by volatility $\sigma$. The Feller condition is violated in (a), while is satisfied in (b).

#### 3.3.3 Computational study

We computationally investigate models with the mean-field component, particularly their randomness and intermittency. We only present the results for the model 2023–2025 because those of the other cases are qualitatively the same. We quantify the difference in randomness inherent to the fish count. The non-negativity-preserving discretization method for the non-negative diffusion bridge [40] directly applies to our case because the average $\bar{X}$ is explicitly known (**Appendix C**). The time interval $(0,1)$ is uniformly divided into 50,000 small time intervals. Each computational case in the sequel uses 1,000,000 sample paths. This numerical method can handle both high- and low-volatility regimes and rigorously preserve the non-negativity of CIR-type processes. If the average is not known, one needs to compute it using another numerical scheme that computes the mean field empirically (e.g., [79, 80]); however, the preservation of nonnegativity of solutions is in general a nontrivial issue, and moreover studies on expectation-dependent (non-Lipschitz) diffusions seem to be scarce. Moreover, existing numerical methods usually assume continuity of drift and diffusion coefficients, but our MVSDE has unbounded coefficients; their ability to preserve the nonnegativity of solutions is unknown [81]. Although the numerical method employed here is specialized for CIR-type processes, it preserves the nonnegativity of solutions and has been computationally found to be convergent [40]. I

**Figure 6(a)** shows the computed probability density functions (PDF) of the MVSDE at several instances, showing that that the PDF is concentrated at the origin. This is due to the intermittency of its solution having a burst-and-rest structure as demonstrated in **Figure 5(a)**. **Figure 6(b)** corresponding to the process of **Figure 5(b)** suggests a unimodal PDF maximized at some positive value, suggesting the dominance of a mean-reverting property. Comparing **Figures 6(a) and 6(b)** shows that the tail of the PDFs are heavier for the case where the diffusion dominates.

Finally, we computationally investigate the behavior of the model when the parameter $\alpha$ controlling the singularity of the diffusion becomes large, i.e., when it does not satisfy **Assumption 1**. Theoretically, the theoretical variance and hence standard deviation of $X$ blows up at the terminal time when $\alpha$ is large as shown in **Figure 7**. This figure also suggests that they are in good agreement when **Assumption 1** is satisfied (**Figures 7(a)–(b)**); however, it is still fluctuating even using 1,000,000 sample paths, which are considered a moderately large number, when **Assumption 1** is satisfied but $\alpha$ is large (**Figure 7(c)**) or is not satisfied (**Figure 7(d)**). Moreover, the computed results without satisfying **Assumption 1** do not present the complete blow-up behavior of the corresponding theoretical ones; instead, the computational standard deviations seem to be approximately vanishing at the terminal time. Increasing the computational resolution and total number of samples may resolve this issue but may be computationally expensive. These findings suggest that theoretical well-posedness study like **Proposition 1** is essential for properly understanding our MVSDE; without it, one may misunderstand that this model admits a pinned solution even for a large $\alpha$ based on Monte-Carlo simulation. Therefore, computational study of SDEs with singularities, such as our MVSDE, should therefore be avoided without theoretical analysis.

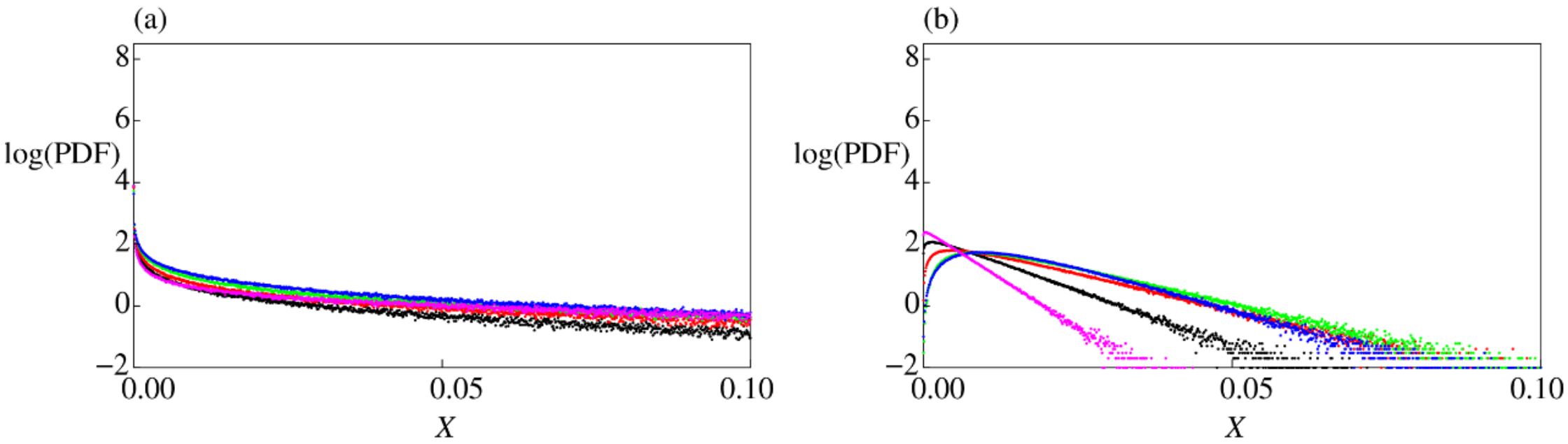

**Figure 6.** Computed PDF of the solution $X$ to the MVSDE for the model of 2023–2025 at several instances, corresponding to **Figures 5(a) and 5(b)**. Colors represent the results at times 0.1 (black), 0.3 (red), 0.5 (green), 0.7 (blue), and 0.9 (magenta). Natural logarithm of PDF is plotted for a visualization purpose.

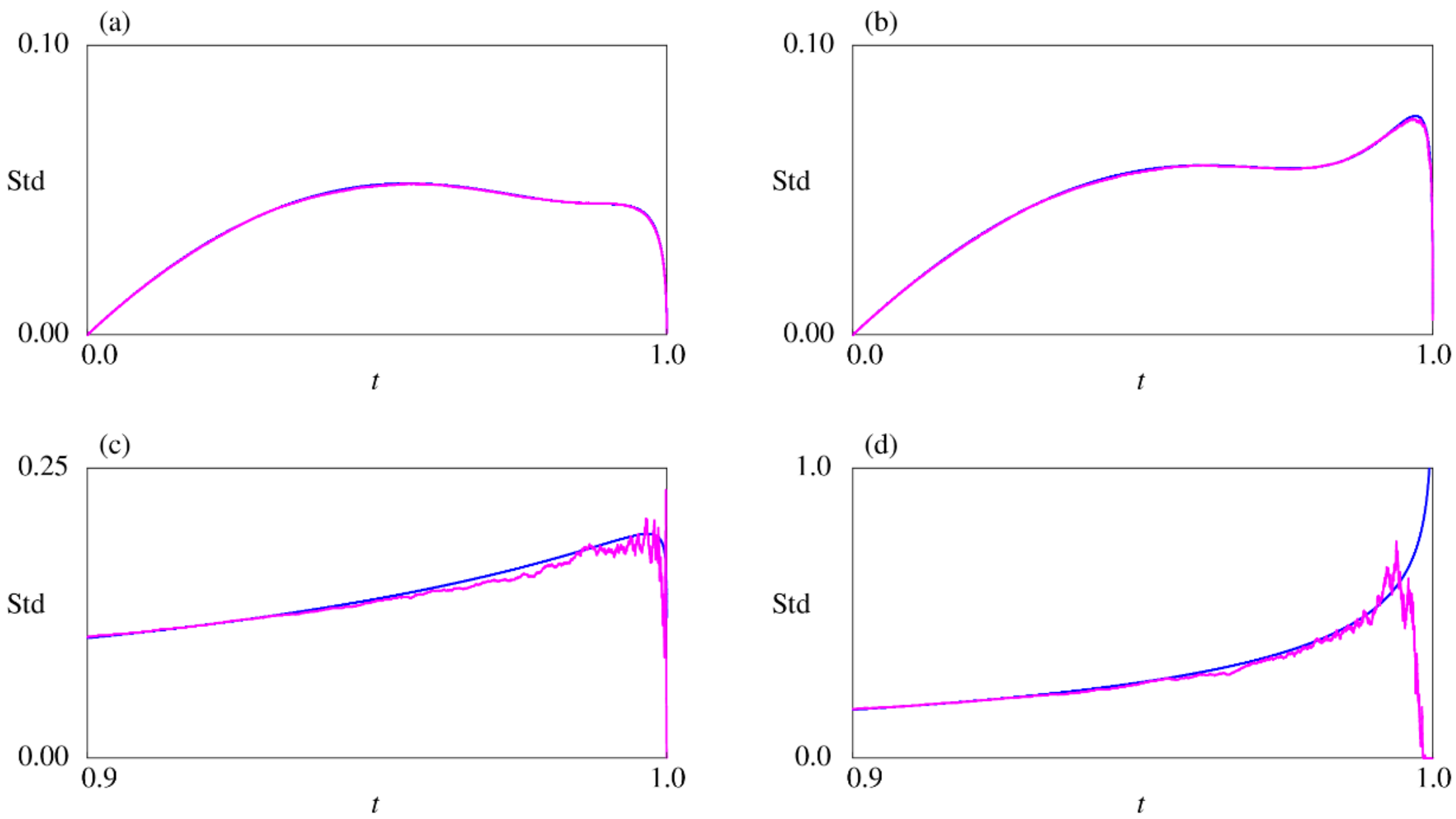

**Figure 7.** Comparison of theoretical and computed standard deviations (Stds) for the model of 2023–2025: the values of $\alpha$ are (a) identified value 0.5482, (b) 1, (c) 1.5, and (d) 2. Only the case (d) does not satisfy **Assumption 1**. Theoretical and computed results are colored blue and magenta, respectively. Only a part of the time interval $(0.9,1)$ is plotted in figure panels (c) and (d) because the focus is around the terminal time.

## 4. Conclusion

We proposed a non-negative diffusion bridge of the MVSDE type, motivated by a pinned CIR process, for modeling the diurnal fish migration phenomenon. We showed that the singularity of the diffusion coefficient must be properly specified to ensure that the MVSDE admits a unique strong solution that is truly pinned at the terminal time. We also discussed the individual-based origin of the MVSDE, showing that this SDE can be interpreted as a superposition of "small" similar SDEs. The MVSDE was applied to the 10-min fish count data of *Ayu* collected from 2023 to 2025 at the study site, and the identified models were found to reasonably fit the empirical data. Applying the Feller condition to the identified models suggested that all models are in the high-volatility regime, where the sample paths are intermittent. We also conducted a computational study to assess the randomness and intermittency of the identified MVSDEs.

Studies on diffusion bridges of MVSDE types are still in an early stage of development. This study investigated a problem that depends only on average, and problems that incorporate higher-order moments or more complex statistics can also be considered at least formally. However, developing such a complex model for practical applications is a prominent issue. For fish migration phenomena, relying on a first-principle method that starts from the space-time swimming behavior of individual fish and then upscales it to an MVSDE may be a possible approach that contributes to both theoretical and applicable studies. Exploring the applications of MVSDE other than fish migration is also an interesting topic. The proposed model for fish migration originates from the CIR process in finance, which connects seemingly disparate objects. We believe several other examples exist worldwide.

This study demonstrated that non-negative diffusion bridges with mean-field effects can be formulated. Because diffusion bridges and related models appear in modern research areas, including machine learning [82] and cell biology [83], we believe that this study opens the possibility of applying mean-field diffusion bridges to various problems not limited to fish migration. Finally, we consider that the superposition approach to derive the MVSDE can be generalized in a way that small processes are heterogeneous, which is currently under investigations by the author.

## Appendix

### A. Proof

***Proof of Proposition 1***

We divide the proof into three steps. The streamline of proof follows those of Section 2 in Bock et al. [45] for mean-field Brownian bridges, but modifications are necessary at each point because the target SDEs are qualitatively different.

***Step 1. The ODE (6)***

We first show that the ODE (6) admits a unique continuous solution (i.e., classical solution that satisfied the ODE pointwise) with $u_0 = u_1 = 0$. As the right-hand side of (6) is Lipschitz continuous with respect to $u_t$ under **Assumption 1**, it admits a unique solution in any interval $[0, 1-\varepsilon)$ for each $\varepsilon \in (0,1)$ (Theorem 6.5 in [84], where we can omit the diffusion term in the proof). Moreover, the coefficient $a$ is non-negative by its definition. Therefore, this solution is non-negative; indeed, if it becomes zero at some time during $t_0 \in (0,1)$ for the first time, then

$$\frac{\mathrm{d}u_t}{\mathrm{d}t} = a(t, u_t) \geq 0 \quad \text{at} \quad t = t_0, \tag{20}$$

indicating that the solution is negative just before $t_0$. This is a contradiction; hence, the solution $u$ is non-negative in any interval $[0, 1-\varepsilon)$ for each $\varepsilon \in (0,1)$.

We present $\lim_{t \to 1-0} u_t = 0$, with which the continuity of the solution on $[0,1]$ is proven. Considering the initial condition $u_0 = 0$, the ODE (6) can be rewritten as follows:

$$u_t = \int_0^t a(s, u_s) \exp\left( -\int_s^t \frac{r(v, u_v)}{1-v} \mathrm{d}v \right) \mathrm{d}s, \quad 0 \leq t < 1. \tag{21}$$

The right-hand side of (21) is evaluated as follows:

$$\begin{aligned} u_t &= \int_0^t a(s, u_s) \exp\left( -\int_s^t \frac{r(v, u_v)}{1-v} \mathrm{d}v \right) \mathrm{d}s \\ &\leq \bar{a} \int_0^t \exp\left( -\int_s^t \frac{r(v, u_v)}{1-v} \mathrm{d}v \right) \mathrm{d}s \\ &\leq \bar{a} \int_0^t \exp\left( -\underline{r} \int_s^t \frac{1}{1-v} \mathrm{d}v \right) \mathrm{d}s \\ &\leq \bar{a} \int_0^t \left( \frac{1-t}{1-s} \right)^{\underline{r}} \mathrm{d}s \end{aligned}, \quad 0 \leq t < 1. \tag{22}$$

An elementary calculation yields

$$\int_0^t \left(\frac{1}{1-s}\right)^{\underline{r}} \mathrm{d}s = \begin{cases} \ln\left(\dfrac{1}{1-t}\right) & (\underline{r}=1) \\ \dfrac{1-(1-t)^{1-\underline{r}}}{1-\underline{r}} & (\underline{r}\neq 1) \end{cases}, \quad 0 \le t < 1 . \tag{23}$$

By (22) and (23), we have

$$\lim_{t\to 1-0} \int_0^t \left(\frac{1-t}{1-s}\right)^{\underline{r}} \mathrm{d}s = 0 . \tag{24}$$

Combining (24) with the non-negativity of the solution $u$ yields the desired result $\lim_{t\to 1-0} u_t = 0$, showing that $u$ is continuous at time 1. This solution $u$ is bounded in $[0,1]$ as well because it is continuous in this domain.

***Step 2. An auxiliary SDE***

We consider the following auxiliary non-negative diffusion bridge, where $u$ denotes a unique solution to the ODE (6):

$$\mathrm{d}Y_t = \left(a(t,u_t) - \frac{r(t,u_t)}{1-t} Y_t\right)\mathrm{d}t + \sigma(t,u_t)\sqrt{\frac{r(t,u_t)}{(1-t)^{\alpha}}|Y_t|}\,\mathrm{d}B_t , \quad 0<t<1 \tag{25}$$

subject to $Y_0 = Y_1 = 0$. The unique existence of strong solutions to the auxiliary SDE (25) in any interval $[0,1-\varepsilon)$ for each $\varepsilon \in (0,1)$ follows from Theorem 6.5 in [84]. Moreover, its non-negativity follows from an argument analogous to that in p. 6-7 in [47] because the coefficients $a(t,u_t), r(t,u_t), \sigma(t,u_t)$ as functions of time $t$ are positive and bounded in $[0,1-\varepsilon)$, and $r$ is strictly bounded away from 0. Using the absolute value of $|Y_t|$ for the diffusion coefficient of (25) thus proves to be superficial. Therefore, in the rest of this proof, we replace "$|Y_t|$" in (25) by "$Y_t$".

Applying a variation of constant formula to (25) yields

$$\begin{aligned} Y_t &= \int_0^t a(s,u_s)\exp\left(-\int_s^t \frac{r(v,u_v)}{1-v}\mathrm{d}v\right)\mathrm{d}s + \int_0^t \sigma(s,u_s)\sqrt{\frac{r(s,u_s)}{(1-s)^{\alpha}}Y_s}\exp\left(-\int_s^t \frac{r(v,u_v)}{1-v}\mathrm{d}v\right)\mathrm{d}B_s , \quad 0 \le t < 1, \\ &\equiv u_t + M_t \end{aligned} \tag{26}$$

where we used (21). The second term $M_t$ in the right-hand side of (26) is a martingale in any interval $[0,1-\varepsilon)$ for each $\varepsilon \in (0,1)$. We show that $\lim_{t\to 1-0} M_t = 0$ with probability 1, with which we obtain $\lim_{t\to 1-0} Y_t = 0$ with probability 1. To show this, it suffices to prove $\lim_{t\to 1-0} \mathbb{E}\left[M_t^2\right] = 0$ because we have $\lim_{t\to 1-0} \mathbb{E}\left[M_t\right] = 0$ from (26). For each $0 \le t \le 1-\varepsilon$, by the isometry (e.g., Proposition 3.20 in [59]), we have

$$
\begin{aligned}
&\mathbb{E}\left[\left(\int_0^t \sigma(s,u_s)\sqrt{\frac{r(s,u_s)}{(1-s)^\alpha}}\,Y_s \exp\left(-\int_s^t \frac{r(v,u_v)}{1-v}\mathrm{d}v\right)\mathrm{d}B_s\right)^2\right] \\
&= \mathbb{E}\left[\int_0^t \sigma^2(s,u_s)\frac{r(s,u_s)}{(1-s)^\alpha}Y_s \exp\left(-2\int_s^t \frac{r(v,u_v)}{1-v}\mathrm{d}v\right)\mathrm{d}s\right] \\
&\le C_\sigma \mathbb{E}\left[\int_0^t \frac{r(s,u_s)}{(1-s)^\alpha}Y_s \exp\left(-2\int_s^t \frac{r(v,u_v)}{1-v}\mathrm{d}v\right)\mathrm{d}s\right] \\
&\le C_\sigma \overline{r}\,\mathbb{E}\left[\int_0^t \frac{1}{(1-s)^\alpha}Y_s \exp\left(-2\int_s^t \frac{r(v,u_v)}{1-v}\mathrm{d}v\right)\mathrm{d}s\right] \\
&\le C_\sigma \overline{r}\,\mathbb{E}\left[\int_0^t \frac{1}{(1-s)^\alpha}Y_s \exp\left(-2\underline{r}\int_s^t \frac{1}{1-v}\mathrm{d}v\right)\mathrm{d}s\right] \\
&= C_\sigma \overline{r}\,\mathbb{E}\left[\int_0^t \frac{1}{(1-s)^\alpha}\left(\frac{1-t}{1-s}\right)^{2\underline{r}} Y_s \mathrm{d}s\right] \\
&= C_\sigma \overline{r}\,(1-t)^{2\underline{r}}\,\mathbb{E}\left[\int_0^t \frac{1}{(1-s)^{\alpha+2\underline{r}}} Y_s \mathrm{d}s\right] \\
&= C_\sigma \overline{r}\,(1-t)^{2\underline{r}}\int_0^t \frac{1}{(1-s)^{\alpha+2\underline{r}}}\mathbb{E}[Y_s]\mathrm{d}s
\end{aligned}
\quad , \tag{27}
$$

where the exchange between expectation and integration in the last line follows from $0 \le t \le 1-\varepsilon$ and integrability of $(Y_t)_{0<t<1}$ (because it is integrable for $0 \le t \le 1-\varepsilon$ due to Theorem 4.4 in [84] with $p=1$ and $\mathbb{E}[Y_t]=u_t$ for $0 \le t < 1$). Here, $C_\sigma = \max_{0\le t\le 1}\sigma^2(t,u_t) \le \max_{0\le s,t\le 1}\sigma^2(s,u_t) < +\infty$ due to the continuity of $\sigma$ and $u$. Because of $\mathbb{E}[Y_t]=u_t$ for $0 \le t < 1$, we can evaluate the last integral in (27) as follows:

$$
\begin{aligned}
\int_0^t \frac{1}{(1-s)^{\alpha+2\underline{r}}}\mathbb{E}[Y_s]\mathrm{d}s &= \int_0^t \frac{1}{(1-s)^{\alpha+2\underline{r}}}u_s \mathrm{d}s \\
&\le \int_0^t \frac{1}{(1-s)^{\alpha+2\underline{r}}}\left(\int_0^s \left(\frac{1-s}{1-v}\right)^{\underline{r}}\mathrm{d}v\right)\mathrm{d}s\,, \\
&= \int_0^t \frac{1}{(1-s)^{\alpha+\underline{r}}}\left(\int_0^s \left(\frac{1}{1-v}\right)^{\underline{r}}\mathrm{d}v\right)\mathrm{d}s
\end{aligned}
\tag{28}
$$

where we used (22). If $\underline{r} \ne 1$, then by (23), we obtain

$$
\int_0^t \frac{1}{(1-s)^{\alpha+2\underline{r}}}\mathbb{E}[Y_s]\mathrm{d}s \le \int_0^t \frac{1}{(1-s)^{\alpha+\underline{r}}}\left(\frac{1-(1-s)^{1-\underline{r}}}{1-\underline{r}}\right)\mathrm{d}s = \frac{1}{1-\underline{r}}\int_0^t \left((1-s)^{-\alpha-\underline{r}} - (1-s)^{1-\alpha-2\underline{r}}\right)\mathrm{d}s\,. \tag{29}
$$

We have

$$
\lim_{t\to 1-0}(1-t)^{2\underline{r}}\int_0^t \left((1-s)^{-\alpha-\underline{r}} - (1-s)^{1-\alpha-2\underline{r}}\right)\mathrm{d}s = 0 \tag{30}
$$

under **Assumption 1** owing to the small singularity condition (5). The case $\underline{r}=1$ follows from the Lôpital's rule. From (27)–(30), we obtain

$$\lim_{t\to 1-0}\mathbb{E}\left[\left(\int_0^t \sigma(s,u_s)\sqrt{\frac{r(s,u_s)}{(1-s)^\alpha}Y_s}\exp\left(-\int_s^t \frac{r(v,u_v)}{1-v}\mathrm{d}v\right)\mathrm{d}B_s\right)^2\right]=0, \quad (31)$$

and hence $\lim_{t\to 1-0}\mathbb{E}\left[M_t^2\right]=0$. Because of $\lim_{t\to 1-0}\mathbb{E}\left[M_t\right]=0$, we obtain $\lim_{t\to 1-0}M_t=0$ and $\lim_{t\to 1-0}Y_t=0$ with probability 1.

Consequently, we can rewrite the SDE (25) as follows

$$\mathrm{d}Y_t=\left(a\left(t,\mathbb{E}\left[Y_t\right]\right)-\frac{r\left(t,\mathbb{E}\left[Y_t\right]\right)}{1-t}Y_t\right)\mathrm{d}t+\sigma\left(t,\mathbb{E}\left[Y_t\right]\right)\sqrt{\frac{r\left(t,\mathbb{E}\left[Y_t\right]\right)}{(1-t)^\alpha}Y_t}\mathrm{d}B_t,\quad 0<t<1. \quad (32)$$

Its unique strong solution satisfies $Y_0=Y_1=0$.

***Step 3. The MVSDE (1)***

The unique solution $Y$ to the auxiliary SDE (25) is a solution to the MVSDE (1) because of representation (32). We prove that this $Y$ is a unique solution to the MVSDE (1). To show this, we check that any solution $X$ to the MVSDE (1) satisfies $\mathbb{E}\left[X_t\right]=u_t$. Taking the expectation in (1) yields

$$\frac{\mathrm{d}}{\mathrm{d}t}\mathbb{E}\left[X_t\right]=a\left(t,\mathbb{E}\left[X_t\right]\right)-\frac{r\left(t,\mathbb{E}\left[X_t\right]\right)}{1-t}\mathbb{E}\left[X_t\right],\quad 0<t<1. \quad (33)$$

From this, we must have $\mathbb{E}\left[X_t\right]=u_t$ because the ODEs (6) and (33) are the same and the solution to the ODE (6) is unique as shown in **Step 1**. Consequently, any strong solution to the MVSDE (1) satisfies $\mathbb{E}\left[X_t\right]=u_t$, hence, it admits a unique strong solution, which is $Y$ in this proof. Its continuity at times zero and one has been proven, and the proof is complete.

□

***Proof of Proposition 2***

First, each SDE (11) admits a unique strong solution with pinned initial and terminal conditions because a version of **Proposition 1** without accounting for the mean-field components applies to it (here, we can consider $\bar{X}$ to be given again according to **Proposition 1**). Then, we can use a calculation analogous to (10). By setting $\hat{X}_t\equiv\sum_{i=1}^{n}X_t^{(i)}$, we have

$$
\begin{aligned}
\mathrm{d}\hat{X}_t &= \sum_{i=1}^{n} \mathrm{d}X_t^{(i)} \\
&= \sum_{i=1}^{n} \left\{ \left( \frac{a\left(t,\bar{X}_t\right)}{n} - \frac{r\left(t,\bar{X}_t\right)}{1-t} X_t^{(i)} \right) \mathrm{d}t + \sigma\left(t,\bar{X}_t\right) \sqrt{\frac{r\left(t,\bar{X}_t\right)}{(1-t)^{\alpha}} X_t^{(i)}} \, \mathrm{d}B_t^{(i)} \right\} \\
&= \left( a\left(t,\bar{X}_t\right) - \frac{r\left(t,\bar{X}_t\right)}{1-t} \sum_{i=1}^{n} X_t^{(i)} \right) \mathrm{d}t + \sigma\left(t,\bar{X}_t\right) \sqrt{\frac{r\left(t,\bar{X}_t\right)}{(1-t)^{\alpha}}} \sum_{i=1}^{n} \sqrt{X_t^{(i)}} \, \mathrm{d}B_t^{(i)} \\
&= \left( a\left(t,\bar{X}_t\right) - \frac{r\left(t,\bar{X}_t\right)}{1-t} \hat{X}_t \right) \mathrm{d}t + \sigma\left(t,\bar{X}_t\right) \sqrt{\frac{r\left(t,\bar{X}_t\right)}{(1-t)^{\alpha}}} \sum_{i=1}^{n} \sqrt{X_t^{(i)}} \, \mathrm{d}B_t^{(i)} \\
&= \left( a\left(t,\bar{X}_t\right) - \frac{r\left(t,\bar{X}_t\right)}{1-t} \hat{X}_t \right) \mathrm{d}t + \sigma\left(t,\bar{X}_t\right) \sqrt{\frac{r\left(t,\bar{X}_t\right)}{(1-t)^{\alpha}}} \sqrt{\sum_{i=1}^{n} X_t^{(i)}} \, \mathrm{d}B_t \\
&= \left( a\left(t,\bar{X}_t\right) - \frac{r\left(t,\bar{X}_t\right)}{1-t} \hat{X}_t \right) \mathrm{d}t + \sigma\left(t,\bar{X}_t\right) \sqrt{\frac{r\left(t,\bar{X}_t\right)}{(1-t)^{\alpha}}} \sqrt{\hat{X}_t} \, \mathrm{d}B_t
\end{aligned}
, \quad 0 < t < 1, \qquad (34)
$$

where $B$ denotes a standard 1-D Brownian motion, and the second last equality is understood in the sense of a law that can be obtained by iteratively applying the procedure employed in (10). The proof is completed because (34) is only the MVSDE (1) with the normalization $T = 1$.

□

## B. Auxiliary data

In **Table B1**, we report fitted parameter values and normalized RMSEs of the theoretical average and standard deviation for models with constant $a, r, \sigma$ with $\alpha = 1$, which are "Model 1" [40].

**Table B1.** Estimated parameter values and normalized RMSEs for the non-mean field case with constant $a, r, \sigma$ and $\alpha = 1$.

| Parameter or index | | 2023 | 2024 | 2025 | 2023-2025 |
|---|---|---|---|---|---|
| $a$ | | 6.003.E-02 | 3.219.E-02 | 3.021.E-02 | 3.673.E-02 |
| $r$ | | 1.647.E+00 | 5.137.E-01 | 4.574.E-01 | 7.100.E-01 |
| $\sigma$ | | 5.942.E-01 | 6.962.E-01 | 7.945.E-01 | 7.252.E-01 |
| Normalized RMSE | Ave | 3.818E-01 | 4.521E-01 | 5.092E-01 | 3.292E-01 |
| | Std | 5.992E-01 | 6.751E-01 | 8.930E-01 | 5.822E-01 |

## C. Theoretical average and variance

The average and variance moments of the unit time fish count for the MVSDE specified in **Section 3** are explicitly given as follows, which can be obtained by taking the expectation of suitable moments of $X$ using (1):

$$
\mathbb{E}\left[X_t\right] = \frac{a}{1-r}\left\{ (1-t)^{r} - (1-t) \right\}, \quad 0 \le t < 1 \qquad (35)
$$

and

$$\mathbb{V}\left[X_t\right]=\frac{ra\mu^2}{1-r}\left\{-\frac{1}{1-r-\alpha}\left((1-t)^{1+r-\alpha}-(1-t)^{2r}\right)+\frac{1}{2-\alpha-2r}\left((1-t)^{2-\alpha}-(1-t)^{2r}\right)\right\}$$
$$+\frac{ra^2\omega}{(1-r)^2}\left\{\begin{array}{l}-\frac{1}{1-\alpha}\left((1-t)^{1+2r-\alpha}-(1-t)^{2r}\right)+\frac{2}{2-r-\alpha}\left((1-t)^{2+r-\alpha}-(1-t)^{2r}\right)\\ -\frac{1}{3-2r-\alpha}\left((1-t)^{3-\alpha}-(1-t)^{2r}\right)\end{array}\right\},\quad 0\leq t<1. \quad (36)$$

Here, we assumed that all the denominators in (35) and (36) do not vanish, but the other case can be obtained by formally taking suitable limits according to Lôpital's rule. Elementary calculations show that $\lim_{t\to 1-0}\mathbb{E}\left[X_t\right]=\lim_{t\to 1-0}\mathbb{V}\left[X_t\right]=0$ under **Assumption 1**.

## References

[1] Foote, K. J., Grant, J. W., & Biron, P. M. Salmonid biomass in streams around the world: a quantitative synthesis. Fish and Fisheries 2025; 26(3): 394–413. https://doi.org/10.1111/faf.12887

[2] Stewart, D. R., Barron, J. C., Harden, T., Grube, E. R., Ulibarri, M., Taylor, A. T., ... & Harris, G. M. The optimal stocking strategy for Yaqui Catfish. North American Journal of Fisheries Management 2023; 43(5): 1407–1426. https://doi.org/10.1002/nafm.10942

[3] Jolly, M. E., Warburton, H. J., Bowie, S., & McIntosh, A. R. Managing isolation: implementing in-stream barriers to exclude introduced trout from fragmented native freshwater fish refuges. River Research and Applications. Online published 2025. https://doi.org/10.1002/rra.4447

[4] Lu, Y., Liu, Q., Li, Y., Jiao, Y., Cheng, B., Sun, G., ... & Qing, J. Measurement and management of fish spawning habitat effectiveness considering potential migration barriers. Ecohydrology 2025; 18(5): e70090. https://doi.org/10.1002/eco.70090

[5] Luiz, O. J., Stratford, D., & Kopf, R. K. Environmental and biological drivers of fish beta diversity and tropical river conservation in northern Australia. Diversity and Distributions 2025; 31(5): e70027. https://doi.org/10.1111/ddi.70027

[6] Colombo, L., & Labrecciosa, P. Resource mobility and market performance. dynamic games and applications 2024; 14(1): 78–96. https://doi.org/10.1007/s13235-023-00517-8

[7] Øksendal, B. Stochastic differential equations: an introduction with applications. Springer 2023, Berlin Heidelberg.

[8] Getz, W. M., Salter, R., Sethi, V., Cain, S., Spiegel, O., & Toledo, S. The statistical building blocks of animal movement simulations. Movement Ecology 2024; 12(1): 67. https://doi.org/10.1186/s40462-024-00507-4

[9] Tao, Y., Giunta, V., Börger, L., & Wilber, M. Q. Towards transient space-use dynamics: re-envisioning models of utilization distribution and their applications. Movement Ecology 2025; 13(1): 12. https://doi.org/10.1186/s40462-025-00538-5

[10] Hanks, E. M., Johnson, D. S., & Hooten, M. B. Reflected stochastic differential equation models for constrained animal movement. Journal of Agricultural, Biological and Environmental Statistics 2017; 22(3): 353–372. https://doi.org/10.1007/s13253-017-0291-8

[11] Pramanik, P. Effects of water currents on fish migration through a Feynman-type path integral approach under 8/3 Liouville-like quantum gravity surfaces. Theory in Biosciences 2021; 140(2): 205–223. https://doi.org/10.1007/s12064-021-00345-7

[12] Yoshioka, H. A stochastic differential game approach toward animal migration. Theory in Biosciences 2019; 138(2): 277–303. https://doi.org/10.1007/s12064-019-00292-4

[13] Yoshioka, H. Superposition of interacting stochastic processes with memory and its application to migrating fish counts. Chaos, Solitons & Fractals 2025-2; 192: 115911. https://doi.org/10.1016/j.chaos.2024.115911

[14] Yoshioka, H., & Yamazaki, K. a jump-driven self-exciting stochastic fish migration model and its fisheries applications. Natural Resource Modeling 2025; 38(1): e12419. https://doi.org/10.1111/nrm.12419

[15] Chu, J., Lam, K. Y., Wang, B., & Wang, T. An optimal switching approach for bird migration. Numerical Algebra, Control and Optimization 2025. Online published. https://doi.org/10.3934/naco.2025026

[16] Cantrell, R. S., Cosner, C., Lam, K. Y., & Mazari-Fouquer, I. Mean field games and ideal free distribution. Journal of Mathematical Biology 2025; 91: 46. https://doi.org/10.1007/s00285-025-02276-z

[17] Knudsen, T. E., MacKenzie, B. R., Thygesen, U. H., & Mariani, P. Evolution and stability of social learning in animal migration. Movement Ecology 2025: 13(1): 43. https://doi.org/10.1186/s40462-025-00564-3

[18] Stella, L., Bauso, D., & Colaneri, P. Mean-field game for collective decision-making in honeybees via switched systems. IEEE Transactions on Automatic Control 2021; 67(8): 3863–3878. https://doi.org/10.1109/TAC.2021.3110166

[19] Borra, F., Cencini, M., & Celani, A. Optimal collision avoidance in swarms of active Brownian particles. Journal of Statistical Mechanics: Theory and Experiment 2021; 2021(8): 083401. https://doi.org/10.1088/1742-5468/ac12c6

[20] Mishura, Y., & Veretennikov, A. Existence and uniqueness theorems for solutions of McKean–Vlasov stochastic equations. Theory of Probability and Mathematical Statistics 2020; 103: 59–101. https://doi.org/10.1090/tpms/1135

[21] Bao, J., & Wang, J. Long-time behavior of one-dimensional McKean-Vlasov SDEs with common noise. Journal of Mathematical Analysis and Applications 2025; 129819. https://doi.org/10.1016/j.jmaa.2025.129819
[22] Pourhadi, E., & Li, C. A uniqueness criterion for McKean–Vlasov fractional stochastic differential equations in $L^p$. Chaos, Solitons & Fractals 2025; 200: 117153. https://doi.org/10.1016/j.chaos.2025.117153
[23] Liu, H., & Lin, J. Y. Stochastic McKean–Vlasov equations with Lévy noise: Existence, attractiveness and stability. Chaos, Solitons & Fractals 2023; 177: 114214. https://doi.org/10.1016/j.chaos.2023.114214
[24] Zhang, H. On a class of Lévy-driven McKean-Vlasov SDEs with Hölder coefficients. Journal of Mathematical Analysis and Applications 2022; 516(2): 126556. https://doi.org/10.1016/j.jmaa.2022.126556
[25] Aïd, R., Bonesini, O., Callegaro, G., & Campi, L. A McKean–Vlasov Game of Commodity Production, Consumption and Trading. Applied Mathematics & Optimization 2022; 86(3): 40. https://doi.org/10.1007/s00245-022-09907-7
[26] Bayer, C., Belomestny, D., Butkovsky, O., & Schoenmakers, J. A reproducing kernel Hilbert space approach to singular local stochastic volatility McKean–Vlasov models. Finance and Stochastics 2024; 28(4): 1147–1178. https://doi.org/10.1007/s00780-024-00541-5
[27] Shrivats, A. V., Firoozi, D., & Jaimungal, S. A mean-field game approach to equilibrium pricing in solar renewable energy certificate markets. Mathematical Finance 2022; 32(3): 779–824.
[28] Hinds, P. D., Sharma, A., & Tretyakov, M. V. Well-posedness and approximation of reflected McKean–Vlasov SDEs with applications. Mathematical Models and Methods in Applied Sciences 2025; 35(08): 1845–1887. https://doi.org/10.1142/S0218202525500241
[29] Huang, H., & Kouhkouh, H. Self-interacting CBO: Existence, uniqueness, and long-time convergence. Applied Mathematics Letters 2025; 161: 109372. https://doi.org/10.1016/j.aml.2024.109372
[30] Lyu, L., & Chen, J. Consensus based stochastic optimal control. In Forty-second International Conference on Machine Learning 2025. https://openreview.net/pdf?id=W7dN3SQKkH
[31] Van Wichelen, J., Verhelst, P., Buysse, D., Belpaire, C., Vlietinck, K., & Coeck, J. Glass eel (Anguilla anguilla L.) behaviour after artificial intake by adjusted tidal barrage management. Estuarine, Coastal and Shelf Science 2021; 249: 107127. https://doi.org/10.1016/j.ecss.2020.107127
[32] Vergeynst, J., Pauwels, I., Baeyens, R., Mouton, A., De Mulder, T., & Nopens, I. Shipping canals on the downstream migration route of European eel (Anguilla anguilla): Opportunity or bottleneck? Ecology of freshwater fish 2021; 30(1): 73–87. https://doi.org/10.1111/eff.12565
[33] Verhelst, P., Westerberg, H., Coeck, J., Harrison, L., Moens, T., Reubens, J., ... & Righton, D. Tidal and circadian patterns of European eel during their spawning migration in the North Sea and the English Channel. Science of the total environment 2023; 905: 167341. https://doi.org/10.1016/j.scitotenv.2023.167341
[34] Harbicht, A. B., Nilsson, P. A., Österling, M., & Calles, O. Environmental and anthropogenic correlates of migratory speeds among Atlantic salmon smolts. River Research and Applications 2021; 37(3): 358–372. https://doi.org/10.1002/rra.3760
[35] Waters, C., Cotter, D., O'Neill, R., Drumm, A., Cooney, J., Bond, N., ... & Maoiléidigh, N. Ó. The use of predator tags to explain reversal movement patterns in Atlantic salmon smolts (Salmo salar L.). Journal of Fish Biology 2025; 106(5): 1316–1333. https://doi.org/10.1111/jfb.15658
[36] Šmejkal, M., Bartoň, D., Blabolil, P., Kolařík, T., Kubečka, J., Sajdlová, Z., ... & Brabec, M. Diverse environmental cues drive the size of reproductive aggregation in a rheophilic fish. Movement Ecology 2023; 11(1): 16. https://doi.org/10.1186/s40462-023-00379-0
[37] Elliott, C. W., Ridgway, M. S., Blanchfield, P. J., & Tufts, B. L. Foraging activity and habitat use throughout an annual migration of adult walleye (Sander vitreus) from the Trent River in eastern Lake Ontario. Animal Biotelemetry 2025; 13(1): 14. https://doi.org/10.1186/s40317-025-00410-8
[38] Ovidio, M., Dierckx, A., & Benitez, J. P. Movement behaviour and fishway performance for endemic and exotic species in a large anthropized river. Limnologica 2023: 99; 126061. https://doi.org/10.1016/j.limno.2023.126061
[39] Tsuji, S., & Shibata, N. Particle size distribution shift and diurnal concentration changes of environmental DNA caused by fish spawning behaviour. Landscape and Ecological Engineering 2025; 21(1): 151–161. https://doi.org/10.1007/s11355-024-00630-9
[40] Yoshioka, H. CIR bridge for modeling of fish migration on sub-hourly scale. Chaos, Solitons & Fractals 2025; 199, Part 3: 116874. https://doi.org/10.1016/j.chaos.2025.116874

[41] Behjoo, H., & Chertkov, M. M. Space-Time Diffusion Bridge. IFAC-PapersOnLine 2024; 58(17): 274–279. https://doi.org/10.1016/j.ifacol.2024.10.181
[42] Chen, Y., & Georgiou, T. Stochastic bridges of linear systems. IEEE Transactions on Automatic Control 2015; 61(2): 526–531. https://doi.org/10.1109/TAC.2015.2440567
[43] Louriki, M. Brownian bridge with random length and pinning point for modelling of financial information. Stochastics 2022; 94(7): 973–1002. https://doi.org/10.1080/17442508.2021.2017438
[44] Sant, J., Jenkins, P. A., Koskela, J., & Spanò, D. EWF: simulating exact paths of the Wright–Fisher diffusion. Bioinformatics 2023; 39(1): btad017. https://doi.org/10.1093/bioinformatics/btad017
[45] Bock, W., Hilbert, A., & Louriki, M. McKean-Vlasov processes of bridge type. Preprint 2025. https://arxiv.org/abs/2501.15568
[46] Cox, J.C., Ingersoll, J. E., & Ross, S. A. A theory of the term structure of interest rates. Econometrica 1985; 53(2): 385–407. https://doi.org/10.2307/1911242
[47] Alfonsi, A. Affine diffusions and related processes: Simulation, theory and applications. Springer, Cham 2015.
[48] Sun, L. H. Systemic risk and interbank lending. Journal of Optimization Theory and Applications 2018, 179(2), 400–424. https://doi.org/10.1007/s10957-017-1185-1
[49] Wang, N., & Zhang, Y. Robust asset-liability management games in a stochastic market with stochastic cash flows under HARA utility. Insurance: Mathematics and Economics 2025, 103125. https://doi.org/10.1016/j.insmatheco.2025.103125
[50] Ham, L., Coomer, M. A., & Stumpf, M. P. The chemical Langevin equation for biochemical systems in dynamic environments. The Journal of Chemical Physics 2022, 157(9). https://doi.org/10.1063/5.0095840
[51] Gorin, G., Vastola, J. J., Fang, M., & Pachter, L. Interpretable and tractable models of transcriptional noise for the rational design of single-molecule quantification experiments. Nature Communications 2022; 13(1): 7620. https://doi.org/10.1038/s41467-022-34857-7
[52] Gorin, G., Vastola, J. J., & Pachter, L. Studying stochastic systems biology of the cell with single-cell genomics data. Cell Systems 2023; 14(10): 822–843. https://doi.org/10.1016/j.cels.2023.08.004
[53] Bossy, M., Jabir, J. F., & Rodriguez, K. M. Instantaneous turbulent kinetic energy modelling based on Lagrangian stochastic approach in CFD and application to wind energy. Journal of Computational Physics 2022; 464: 110929. https://doi.org/10.1016/j.jcp.2021.110929
[54] Yoshioka, H., & Yoshioka, Y. Non-Markovian superposition process model for stochastically describing concentration–discharge relationship. Chaos, Solitons & Fractals 2025; 199: 116715. https://doi.org/10.1016/j.chaos.2025.116715
[55] Bihun, C. J., Faust, M. D., Kraus, R. T., MacDougall, T. M., Robinson, J. M., Vandergoot, C. S., & Raby, G. D. Is sexual size dimorphism in walleye, Sander vitreus, a driver of seasonal movements in Lake Erie?. Journal of Fish Biology 2025; 106(2): 430–441. https://doi.org/10.1111/jfb.15960
[56] Fujita, M. H., Binder, T. R., Henderson, M., & Marsden, J. E. Modeling regional occupancy of fishes using acoustic telemetry: a model comparison framework applied to lake trout. Animal Biotelemetry 2024; 12(1): 25. https://doi.org/10.1186/s40317-024-00380-3
[57] Swadling, D. S., Knott, N. A., Taylor, M. D., Rees, M. J., Cadiou, G., & Davis, A. R. Consequences of Juvenile fish movement and seascape connectivity: Does the concept of nursery habitat need a rethink?. Estuaries and Coasts 2024; 47(3): 607–621. https://doi.org/10.1007/s12237-023-01323-6
[58] Sato, D., & Seguchi, Y. Estimation of Ayu-fish migration status and passage time through the Yodo river barrage fishway. Advances in River Engineering 2024; 30: 1–4. In Japanese with English Abstract. https://doi.org/10.11532/river.30.0_1
[59] Capasso, V., & Bakstein, D. An introduction to continuous-time stochastic processes: Theory, models, and applications to finance, biology, and medicine. Springer, Cham 2021.
[60] Bjørnås, K. L., Railsback, S. F., Calles, O., & Piccolo, J. J. Modeling Atlantic salmon (Salmo salar) and brown trout (S. trutta) population responses and interactions under increased minimum flow in a regulated river. Ecological Engineering 2021; 162: 106182. https://doi.org/10.1016/j.ecoleng.2021.106182
[61] Gao, D., Zhu, X., Huang, M., Wang, S., & Guo, H. Modeling fish swimming trajectories in a sudden expansion flow based on Eulerian Lagrangian agent method (ELAM): A case study of red crucian carp. Ecological Modelling 2025; 510: 111286. https://doi.org/10.1016/j.ecolmodel.2025.111286
[62] Ruiz-Coello, M. X., Bottacin-Busolin, A., & Marion, A. A Eulerian-Lagrangian model for simulating fish pathline distributions in vertical slot fishways. Ecological Engineering 2024; 203: 107264. https://doi.org/10.1016/j.ecoleng.2024.107264

[63] Hajiesmaeili, M. Scaling up individual-based models for salmonid populations in hydropower-regulated rivers: key advances and challenges. River Research and Applications 2025. Online published. https://doi.org/10.1002/rra.70023
[64] Dahl, K. R., & Eyjolfsson, H. Self-exciting jump processes and their asymptotic behaviour. Stochastics 2022; 94(8): 1166–1185. https://doi.org/10.1080/17442508.2022.2028789
[65] Schmutz, V., Löcherbach, E., & Schwalger, T. On a finite-size neuronal population equation. SIAM Journal on Applied Dynamical Systems 2023; 22(2): 996–1029. https://doi.org/10.1137/21M1445041
[66] Ida, S., Sato, T., Sueyoshi, M., Kishi D., & Ohta, T. Construction of an isoscape of the strontium isotope ratios of the entire Nagara River in central Japan. Limnology 2025; 26: 549–556. https://doi.org/10.1007/s10201-025-00787-8
[67] Tsukamoto, K., & Uchida, K. (1992). Migration mechanism of the ayu. In Oceanic and Anthropogenic Controls of Life in the Pacific Ocean: Proceedings of the 2nd Pacific Symposium on Marine Sciences 1992, Nadhodka, Russia, August 11–19, 1988 (pp. 145–172). Springer Netherlands, Dordrecht. https://doi.org/10.1007/978-94-011-2773-8_12
[68] Jadhav, V., Guttal, V., & Masila, D. R. Randomness in the choice of neighbours promotes cohesion in mobile animal groups. Royal Society Open Science 2022; 9(3): 220124. https://doi.org/10.1098/rsos.220124
[69] Shea, J., & Stark, H. Emergent collective behavior of cohesive, aligning particles. The European Physical Journal E 2025; 48(4): 22. https://doi.org/10.1140/epje/s10189-025-00482-7
[70] Monk, C. T., Aslak, U., Brockmann, D., & Arlinghaus, R. Rhythm of relationships in a social fish over the course of a full year in the wild. Movement Ecology 2023; 11(1): 56. https://doi.org/10.1186/s40462-023-00410-4
[71] Wang, W., Escobedo, R., Sanchez, S., Han, Z., Sire, C., & Theraulaz, G. Collective phases and long-term dynamics in a fish school model with burst-and-coast swimming. Royal Society Open Science 2025; 12(5): 240885. https://doi.org/10.1098/rsos.240885
[72] Bicego, S., Kalise, D., & Pavliotis, G. A. Computation and control of unstable steady states for mean field multiagent systems. The Royal Society Proceedings A 2025; 481(2311): 20240476. https://doi.org/10.1098/rspa.2024.0476
[73] Karatzas, I., & Shreve, S. Brownian motion and stochastic calculus. Springer, New York.
[74] Bučková, Z., Ehrhardt, M., & Günther, M. Fichera theory and its application in finance. In European Consortium for Mathematics in Industry 2014; pp. 103–111. Springer, Cham. https://doi.org/10.1007/978-3-319-23413-7_13
[75] Ceci, C., Bufalo, M., & Orlando, G. Modelling the industrial production of electric and gas utilities through the $CIR^3$ model. Mathematics and Financial Economics 2023; 18: 1–25. https://doi.org/10.1007/s11579-023-00350-y
[76] Pagliarani, S., & Pascucci, A. The exact Taylor formula of the implied volatility. Finance and Stochastics 2017; 21(3): 661–718. https://doi.org/10.1007/s00780-017-0330-x
[77] Veraart, A. E., & Veraart, L. A. Stochastic volatility and stochastic leverage. Annals of Finance 2012, 8(2), 205–233. https://doi.org/10.1007/s10436-010-0157-3
[78] Xu, G., & Wang, Y. On stability of the Markov-modulated skew CIR process. Statistics & Probability Letters 2016;109: 139–144. https://doi.org/10.1016/j.spl.2015.10.020
[79] Nobis, G., Belova, A., Springenberg, M., Daems, R., Knochenhauer, C., Opper, M., ... & Samek, W. Fractional Brownian Bridges for Aligned Data 2025, March. In Learning Meaningful Representations of Life (LMRL) Workshop at ICLR 2025. https://openreview.net/pdf?id=PoEFXbYNii
[80] Le Treut, G., Ancheta, S., Huber, G., Orland, H., & Yllanes, D. Markov-bridge generation of transition paths and its application to cell-fate choice. Physical Review Research 2025; 7(1): 013010. https://doi.org/10.1103/PhysRevResearch.7.013010
[81] Jie, L., Luo, L., & Zhang, H. One-dimensional McKean–Vlasov stochastic Volterra equations with Hölder diffusion coefficients. Statistics & Probability Letters 2024; 205: 109970. https://doi.org/10.1016/j.spl.2023.109970
[82] Liu, H., Shi, B., & Wu, F. Tamed Euler–Maruyama approximation of McKean–Vlasov stochastic differential equations with super-linear drift and Hölder diffusion coefficients. Applied Numerical Mathematics 2023; 183: 56-85. https://doi.org/10.1016/j.apnum.2022.08.012
[83] Zhang, J., Wu, D., Li, Z., & Xu, L. A class of time-changed McKean-Vlasov stochastic differential equations with super-linear drift and Hölder diffusion coefficients. Communications in Nonlinear Science and Numerical Simulation 2025; 152: 109304. https://doi.org/10.1016/j.cnsns.2025.109304
[84] Mao, X. (2007). Stochastic differential equations and applications. Elsevier.